\documentclass[12pt,a4paper]{article}%
\usepackage{amsmath}
\usepackage{amsthm}
\usepackage{amssymb}
\usepackage{amsfonts}
\usepackage[{left=1.5cm, right=1.5cm, top=3cm, bottom=3cm}]{geometry}
\usepackage{graphicx}%
\providecommand{\U}[1]{\protect\rule{.1in}{.1in}}
\newtheorem{theorem}{Theorem}

\newtheorem{corollary}{Corollary}[theorem]

\newtheorem{lemma}{Lemma}

\newtheorem{proposition}{Proposition}

\theoremstyle{definition}

\newtheorem{example}{Example}

\begin{document}

\title{Super Catalan Numbers and \\Fourier Summation over Finite Fields}
\author{Kevin Limanta\\School of Mathematics and Statistics\\UNSW Sydney\\\texttt{k.limanta@unsw.edu.au}
\and Norman Wildberger\\School of Mathematics and Statistics\\UNSW Sydney\\\texttt{n.wildberger@unsw.edu.au}}
\date{}
\maketitle

\begin{abstract}
We study polynomial summation over unit circles over finite fields of odd
characteristic, obtaining a purely algebraic integration theory without
recourse to infinite procedures. There are nonetheless strong parallels to
classical integration theory over a circle, and we show that the super Catalan
numbers and closely related rational numbers lie at the heart of both
theories. This gives a uniform analytic meaning to these up to now somewhat
mysterious numbers. 

Our derivation utilises the three-fold symmetry of chromogeometry between
Euclidean and relativistic geometries, and we find that the Fourier summation
formulas we derive in these two different settings are closely connected.


\end{abstract}

\section{Introduction}

This paper will elucidate the crucial and perhaps surprising role that super
Catalan numbers, and closely associated rational numbers which we will call
circular super Catalan numbers, play in algebra and analysis. \textbf{Super
Catalan numbers} are natural numbers of the form
\begin{equation}
S(m,n)\equiv\frac{\left(  2m\right)  !\left(  2n\right)  !}{m!n!\left(
m+n\right)  !} \nonumber%
\end{equation}
and they generalize the Catalan numbers, since $S\left(  1,n\right)  $ is
twice the Catalan number $C\left(  n\right)  $. While the super Catalan
numbers were first introduced by Catalan \cite{Catalan} in 1874, it appears
that the first modern study of these numbers was initiated by Gessel
\cite{Gessel} in 1992.

Clearly $S\left(  m,n\right)  $ is symmetric in $m$ and $n$, and $S\left(
0,n\right)  =S\left(  n,n\right)  =\binom{2n}{n}$. The fact that $S(m,n)$ is
always an integer is a simple consequence of the recurrence relation
$4S(m,n)=S(m+1,n)+S(m,n+1)$. The table below lists $S(m,n)$ for $0\leq
m,n\leq10$.%

\[%
\begin{tabular}
[c]{c|ccccccccccc}%
$m\backslash n$ & $0$ & $1$ & $2$ & $3$ & $4$ & $5$ & $6$ & $7$ & $8$ & $9$ &
$10$\\\hline
$0$ & $1$ & $2$ & $6$ & $20$ & $70$ & $252$ & $924$ & $3432$ & $12\,870$ &
$48\,620$ & $184\,756$\\
$1$ & $2$ & $2$ & $4$ & $10$ & $28$ & $84$ & $264$ & $858$ & $2860$ &
$\allowbreak9724$ & $33\,592$\\
$2$ & $6$ & $4$ & $6$ & $12$ & $28$ & $72$ & $198$ & $572$ & $1716$ & $5304$ &
$\allowbreak16\,796$\\
$3$ & $20$ & $10$ & $12$ & $20$ & $40$ & $90$ & $220$ & $572$ & $1560$ &
$4420$ & $12\,920$\\
$4$ & $70$ & $28$ & $28$ & $40$ & $70$ & $140$ & $308$ & $728$ & $1820$ &
$4760$ & $\allowbreak12\,920$\\
$5$ & $252$ & $84$ & $72$ & $90$ & $140$ & $252$ & $504$ & $1092$ & $2520$ &
$6120$ & $15\,504$\\
$6$ & $924$ & $264$ & $198$ & $220$ & $308$ & $504$ & $924$ & $1848$ & $3960$
& $8976$ & $21\,318$\\
$7$ & $3432$ & $858$ & $572$ & $572$ & $728$ & $1092$ & $1848$ & $3432$ &
$6864$ & $14\,586$ & $32\,604$\\
$8$ & $12\,870$ & $2860$ & $1716$ & $1560$ & $1820$ & $2520$ & $3960$ & $6864$
& $12\,870$ & $25\,740$ & $54\,340$\\
$9$ & $48\,620$ & $\allowbreak9724$ & $\allowbreak5304$ & $4420$ & $4760$ &
$\allowbreak6120$ & $8976$ & $14\,586$ & $25\,740$ & $48\,620$ &
$\allowbreak97\,240$\\
$10$ & $184\,756$ & $33\,592$ & $\allowbreak16\,796$ & $12\,920$ &
$\allowbreak12\,920$ & $15\,504$ & $21\,318$ & $32\,604$ & $54\,340$ &
$\allowbreak97\,240$ & $184\,756$%
\end{tabular}
\ \
\]

\begin{center}
\textbf{Table 1}: $S(m,n)$ for $0\leq m,n\leq10$.
\end{center}

Somewhat surprisingly, there is no known general combinatorial understanding
for $S\left(  m,n\right)  $ to date, unlike the Catalan numbers (sequence
A000108 in the OEIS \cite{OEIS}) which have hundreds of combinatorial
interpretations, compiled for example in \cite{Stanley}. There are however
some combinatorial or weighted interpretations of $S\left(  m,n\right)  $ for
specific values of $m$ and $n;$ in terms of blossom trees \cite{Schaeffer},
cubic trees \cite{Pippenger}, or pairs of Dyck paths with restricted heights
for $\left(  m,n\right)  =\left(  2,3\right)  $ by Gessel \cite{Gessel2}. 

Chen and Wang gave an interpretation for $S(m,m+s)$ for $0\leq s\leq3$ in
terms of restricted lattice paths \cite{Chen}. Recently, a combinatorial
interpretation for $S(3,n)$ and $S(4,n)$ was given by Gheorghiciuc and
Orelowitz \cite{Gheorghiciuc} as the numbers of $3$-tuples and $4$-tuples of
Dyck paths satisfying certain conditions.

A weighted interpretation of $S(m,n)$ involve counting $2$-Motzkin paths was
discovered by Allen and Gheorghiciuc \cite{Allen}, and another in terms of
certain values of Krawtchouk polynomials was found by Georgiadis, Munemasa,
and Tanaka in \cite{Georgiadis}. The latter will play a prominent role in this paper.

Finally, there are some other established identities involving the super
Catalan numbers, for example relations to the reciprocal super Catalan matrix
were discovered by Prodinger \cite{Prodinger}. Liu \cite{Liu} found an
identity for congruences of the super Catalan numbers, and the paper by
Miki\'{c} \cite{Mikic} gives an identity for convolution of the super Catalan numbers.

In this paper we introduce also a closely related family of rational numbers,
the \textbf{circular super Catalan numbers }given by
\begin{equation}
\Omega(m,n)\equiv\frac{S(m,n)}{4^{m+n}}
\nonumber%
\end{equation}
and a table of $\Omega(m,n)$ for small values of $m$ and $n$ is given below.%
\[%
\begin{tabular}
[c]{c|ccccccccc}%
$m\backslash n$ & $0$ & $1$ & $2$ & $3$ & $4$ & $5$ & $6$ & $7$ & $8$\\\hline
$0$ & $1$ & $\frac{1}{2}$ & $\frac{3}{8}$ & $\frac{5}{16}$ & $\frac{35}{128}$
& $\frac{63}{256}$ & $\frac{231}{1024}$ & $\frac{429}{2048}$ & $\frac
{6435}{32\,768}$\\
$1$ & $\frac{1}{2}$ & $\frac{1}{8}$ & $\frac{1}{16}$ & $\frac{5}{128}$ &
$\frac{7}{256}$ & $\frac{21}{1024}$ & $\frac{33}{2048}$ & $\frac{429}%
{32\,768}$ & $\frac{715}{65\,536}$\\
$2$ & $\frac{3}{8}$ & $\frac{1}{16}$ & $\frac{3}{128}$ & $\frac{3}{256}$ &
$\frac{7}{1024}$ & $\frac{9}{2048}$ & $\frac{99}{32\,768}$ & $\frac
{143}{65\,536}$ & $\frac{429}{262\,144}$\\
$3$ & $\frac{5}{16}$ & $\frac{5}{128}$ & $\frac{3}{256}$ & $\frac{5}{1024}$ &
$\frac{5}{2048}$ & $\frac{45}{32\,768}$ & $\frac{55}{65\,536}$ & $\frac
{143}{262\,144}$ & $\frac{195}{524\,288}$\\
$4$ & $\frac{35}{128}$ & $\frac{7}{256}$ & $\frac{7}{1024}$ & $\frac{5}{2048}$
& $\frac{35}{32\,768}$ & $\frac{35}{65\,536}$ & $\frac{77}{262\,144}$ &
$\frac{91}{524\,288}$ & $\frac{455}{4194\,304}$\\
$5$ & $\frac{63}{256}$ & $\frac{21}{1024}$ & $\frac{9}{2048}$ & $\frac
{45}{32\,768}$ & $\frac{35}{65\,536}$ & $\frac{63}{262\,144}$ & $\frac
{63}{524\,288}$ & $\frac{273}{4194\,304}$ & $\frac{315}{8388\,608}$\\
$6$ & $\frac{231}{1024}$ & $\frac{33}{2048}$ & $\frac{99}{32\,768}$ &
$\frac{55}{65\,536}$ & $\frac{77}{262\,144}$ & $\frac{63}{524\,288}$ &
$\frac{231}{4194\,304}$ & $\frac{231}{8388\,608}$ & $\frac{495}{33\,554\,432}%
$\\
$7$ & $\frac{429}{2048}$ & $\frac{429}{32\,768}$ & $\frac{143}{65\,536}$ &
$\frac{143}{262\,144}$ & $\frac{91}{524\,288}$ & $\frac{273}{4194\,304}$ &
$\frac{231}{8388\,608}$ & $\frac{429}{33\,554\,432}$ & $\frac{429}%
{67\,108\,864}$\\
$8$ & $\frac{6435}{32\,768}$ & $\frac{715}{65\,536}$ & $\frac{429}{262\,144}$
& $\frac{195}{524\,288}$ & $\frac{455}{4194\,304}$ & $\frac{315}{8388\,608}$ &
$\frac{495}{33\,554\,432}$ & $\frac{429}{67\,108\,864}$ & $\frac
{6435}{2147\,483\,648}$%
\end{tabular}
\ \ \
\]

\begin{center}
\textbf{Table 2}: $\Omega(m,n)$ for $0\leq m,n\leq8$.
\end{center}

We see that the symmetry property for $\Omega(m,n)$ is still retained, and the
recurrence property of the super Catalan numbers becomes the Pascal-like
property
\[
\Omega(m,n)=\Omega(m+1,n)+\Omega(m,n+1).
\]
It is not hard to see that the table above is completely determined by the
values in the zeroth row, the symmetry property, and the Pascal-like property above.

The objective of this paper is to show that the circular super Catalan numbers
have a central algebraic and analytic interpretation in the context of
integration, or summation, theory of polynomials over unit circles in finite
fields of odd characteristic. This role extends even to summations
over suitably defined unit circles in relativistic geometries, and the curious
relations between the Euclidean and relativistic formulas are an aspect of the
three-fold symmetry of \textit{chromogeometry} \cite{Wildberger}. 

This theory brings together Euclidean planar geometry with quadratic form
$x^{2}+y^{2}$ (blue) with two different forms of relativistic
Einstein-Minkowski geometries with quadratic forms $x^{2}-y^{2}$ (red) and
$xy$ (green). The relations and symmetries between the three geometries play a
crucial role in our derivations, as summations over unit circles in both
Euclidean and relativistic geometries are closely related, and connections
between them provide a powerful tool for explicit evaluations.

We briefly summarize the remarkable formulas in which $\Omega(m,n)$ play the
primary role and illustrate them with examples below. We will also observe
that the classical formulas for integration on the circle in the Euclidean
plane can, somewhat surprisingly, also be restated in terms of super Catalan
numbers, so that we obtain a bridge between characteristic zero and prime
characteristic harmonic analysis. 

We will see that these formulas also suggest the study of families of
additional rational numbers that play secondary roles.

\subsection{Summing polynomials over unit circles over finite fields}

The problem of summing polynomials over a finite field has many implications
to number theory, Fourier analysis, and more generally harmonic analysis. In
2015, Thakur gave an expository survey on sums of integral powers of monic
polynomials of certain degrees over finite fields \cite{Thakur}, in which
relations to the zeta function also appear. For work aimed more in the
direction of complex number values see \cite{Bourgain}, \cite{Carbery}, \cite{Gurevich}, and \cite{Tao}.

To briefly introduce the key ideas in our approach, we begin with a finite
field $\mathbb{F}_{q}$ with multiplicative identity $1_{q},$ where $q=p^{r}$
is a power of an odd prime $p,$ so that $\mathbb{F}_{p}\subseteq\mathbb{F}%
_{q}.$ The associated affine plane $\mathbb{A=A}\left(  \mathbb{F}_{q}\right)
$ consists of points $\left[  x,y\right]  $ where $x,y\in\mathbb{F}_{q}$. We
introduce three different metrical geometries: the Euclidean and two
relativistic geometries which we color-code as blue ($b$), red ($r$), and green
($g$) with respective associated quadratic forms $x^{2}+y^{2}$, $x^{2}-y^{2}$,
and $xy,$ and where the subscript $c$ refers to any of these. The unit circles
in each geometry are then defined by
\begin{align*}
S_{b}\left(\mathbb{F}_{q}\right) & \equiv\left\{\left[x,y\right]
\in\mathbb{A}\ \colon\ x^{2}+y^{2}=1\right\}  \\ S_{r}\left(\mathbb{F}_{q}\right)  &\equiv\left\{  \left[x,y\right]  \in\mathbb{A}
\ \colon\ x^{2}-y^{2}=1\right\}  \\
S_{g}\left(\mathbb{F}_{q}\right) & \equiv\left\{\left[x,y\right]
\in\mathbb{A}\ \colon\ xy=1\right\}.
\end{align*}
We will use $S_{c,q}$ as a shorthand for $S_{c}\left(\mathbb{F}_{q}\right)
$ when the underlying field is clear from context.

Recall that for a prime power $q=p^{r}$ the \textbf{Jacobi symbol} is defined
as
\begin{align*}
\left(\frac{-1}{q}\right)  \equiv\left(\frac{-1}{p}\right)^{r}%
\end{align*}
where the notation $\left(\frac{-1}{p}\right)$ refers to the Legendre
symbol, which has value $1$ if $-1$ is a quadratic residue $\operatorname{mod}%
p$, that is a square in the field $\mathbb{F}_{p}$, and has value of $-1$ if
$-1$ is not a quadratic residue. The Jacobi symbol comes into play in the
following elementary count of the sizes of the unit circles, which we present in the order green, red, and blue here to align with the natural
increasing complexity of formulas.

\begin{lemma}
\label{Lemma size of unit circles Intro}Let $S_{b,q}$, $S_{r,q}$, and
$S_{g,q}$ be the unit circles in the blue, red, and green geometries
respectively over the field $\mathbb{F}_{q}$. Then%
\[
\left\vert S_{g,q}\right\vert =q-1\qquad\left\vert S_{r,q}\right\vert
=q-1\quad\mathrm{and}\quad\left\vert S_{b,q}\right\vert =q-\left(  \frac
{-1}{q}\right).
\]
\end{lemma}

Polynumbers are a variant of polynomials introduced by Wildberger in his online \textit{Maths Foundation} lectures, and can be defined without having a prior notion of
variables. This allows us to more clearly distinguish between \textit{formal
polynomials} and \textit{polynomial functions}, and allows us to build an integration
or summation theory over finite fields in a cleaner and more rigorous way. We
may think of a polynumber as a list or array of coefficients, equipped with
the usual Cauchy product, giving us the algebra $\mathrm{Pol}_{2}\left(
\mathbb{F}_{q}\right)  $. We can evaluate a polynumber $\pi$ to get a
polynomial function $\varepsilon\left(  \pi\right)  $, so there is an
evaluation map $\pi\mapsto\varepsilon\left(  \pi\right)  $.

For each geometry we define a linear functional on $\mathrm{Pol}_{2}\left(
\mathbb{F}_{q}\right)  ,$ called the \textbf{Fourier summation functional}, as
follows:
\begin{align}
\psi_{c,q} \enspace \colon \enspace \mathrm{Pol}_{2}\left(\mathbb{F}_{q}\right) &\rightarrow
\mathbb{F}_{q}\nonumber\\
\pi &\mapsto\frac{1}{\left\vert S_{c,q}\right\vert }\sum_{\left[
x,y\right]  \in S_{c,q}}\varepsilon\left(  \pi\right)  \left(  x,y\right)
\nonumber%
\end{align}
where we treat $\left\vert S_{c,q}\right\vert $ as an element of
$\mathbb{F}_{q}$. Its reciprocal is well-defined since $\left\vert
S_{c,q}\right\vert \neq0$ from Lemma \ref{Lemma size of unit circles Intro}.
We will show that $\psi_{c,q}$ satisfies the following three conditions:

\begin{description}
\item[(Normalization)] For the multiplicative identity $\mathbf{1}%
\in\mathrm{Pol}_{2}\left(  \mathbb{F}_{q}\right)  $, we have $\psi
_{c,q}\left(  \mathbf{1}\right)  =1_{q}$.

\item[(Locality)] If $\pi$ is a polynumber for which $\varepsilon\left(
\pi\right)  $ is the zero function when restricted to $S_{c,q}$, then
$\psi_{c,q}\left(  \pi\right)  =0$.

\item[(Invariance)] $\psi_{c,q}$ is rotationally invariant with respect to
the group of linear isometries of the affine space determined by the geometry.
\end{description}

Any linear functional satisfying the three conditions above will be called a
\textbf{circular integral functional}. A main result is the following
uniqueness theorem.

\begin{theorem}
\label{Thm Uniqueness Intro}The Fourier summation functional $\psi_{c,q}$ is
the unique circular integral functional in the geometry with color $c$, for
each of the green, red, and blue geometries.
\end{theorem}

The second half of this paper is dedicated to finding an explicit formula for
$\psi_{c,q}$ in each case. The result in the blue (Euclidean) geometry is
established from those in the green and red (relativistic) geometries. The
explicit formulas for $\psi_{c,q}$ in the red and blue geometries involve the
circular super Catalan numbers $\Omega(m,n)$ reduced modulo $p$, which we
think of as elements of $\mathbb{F}_{p}$. Note that the expression
\[
\Omega(m,n)\operatorname{mod}p=\frac{S\left(  m,n\right)  }{4^{m+n}%
}\operatorname{mod}p
\]
is well-defined since $S\left(  m,n\right)  $ is a natural number and
$4^{m+n}$ is never a multiple of $p$. A simplified and restricted version of
the formulas for $\psi_{c,q}$ is given below.

\begin{theorem}
\label{Thm CIF formula SC Intro}Consider a finite field $\mathbb{F}_{q}$ where
$q=p^{r}$ for some prime $p>2$.

\begin{enumerate}
\item In the red geometry, for $0\leq k+l<q-1$,
\[
\psi_{r,q}\left(  \alpha^{k}\beta^{l}\right)  =\left\{
\begin{array}
[c]{cl}%
\left(  -1\right)  ^{n}\Omega(m,n)\operatorname{mod}p & \text{if }k=2m\text{
and }l=2n\text{,}\\
0 & \text{otherwise.}%
\end{array}
\right.
\]

\item In the blue geometry, for $0\leq k+l<q-1$,
\[
\psi_{b,q}\left(  \alpha^{k}\beta^{l}\right)  =\left\{
\begin{array}
[c]{cl}%
\Omega(m,n)\operatorname{mod}p & \text{if }k=2m\text{ and }l=2n\text{,}\\
0 & \text{otherwise.}%
\end{array}
\right.
\]

\end{enumerate}
\end{theorem}

As a consequence, the normalized sums $\psi_{r,q}\left(  \alpha^{k}\beta
^{l}\right)  $ and $\psi_{b,q}\left(  \alpha^{k}\beta^{l}\right)  $ lie in the
prime field $\mathbb{F}_{p},$ and are independent of the power $q$ of $p$.
Moreover, $\psi_{r,q}\left(  \alpha^{2m}\beta^{2n}\right)  $ and $\psi
_{b,q}\left(  \alpha^{2m}\beta^{2n}\right)  $ are never $0$. Perhaps most
remarkably, the single rational number $\Omega(m,n)$ determines the values
$\psi_{r,q}\left(  \alpha^{2m}\beta^{2n}\right)  $ and $\psi_{b,q}\left(
\alpha^{2m}\beta^{2n}\right)  $ for \textit{all} prime powers $q>2\left(
m+n\right)  +1$.

Here is an example to illustrate these formulas.
\begin{example}
\label{Example q = 13} We work with $\mathbb{F}_{13}$ which we identify as
$\left\{0,1,2,\ldots,12\right\}$. Here, $-1 = 12$ is a square so by Lemma \ref{Lemma size of unit circles Intro}, $\left\vert S_{b,13}\right\vert =12=-1$ regarded as an element of
$\mathbb{F}_{13}$. Using Theorem \ref{Thm CIF formula SC Intro}, if $0\leq
m+n<6$ then $\psi_{b,13}\left(\alpha^{2m}\beta^{2n}\right) = \Omega(m,n)\operatorname{mod}13$. The blue unit circle is 
\begin{align*}
    S_{b,13}=\left\{
\left[  1,0\right]  \ \left[  12,0\right]  \ \left[  0,1\right]  \ \left[
0,12\right]  \ \left[  2,6\right]  \ \left[  2,7\right]  \ \left[
11,6\right]  \ \left[  11,7\right]  \ \left[  6,2\right]  \ \left[
7,2\right]  \ \left[  6,11\right]  \ \left[  7,11\right]  \right\}
\end{align*}
where there is no need for commas as the points are not naturally ordered. By manual
inspection, and omitting the zero terms,
\begin{align*}
\psi_{b,13}\left(  \alpha^{2}\beta^{6}\right)   &  =-\sum_{\left[  x,y\right]
\in S_{b,13}}x^{2}y^{6}\\
&  =-\left(  2^{2}6^{6}+2^{2}7^{6}+11^{2}6^{6}+11^{2}7^{6}+6^{2}2^{6}%
+7^{2}2^{6}+6^{2}11^{6}+7^{2}11^{6}\right)  =4.
\end{align*}
We may verify that $\Omega(1,3)\operatorname{mod}13=\frac{5}{128}%
\operatorname{mod}13=4$.
\end{example}

The complete formulas for the red and blue geometries which go beyond the
constraint $k+l < q-1$ are more involved, and they point towards an additional
family of integers that extend the super Catalan numbers, and that arise as
coefficients of what we call the \textbf{circular polynumbers}
\[
\pi_{k,l}\equiv\left(  \frac{1+\alpha}{2}\right)  ^{k}\left(  \frac{1-\alpha
}{2}\right)  ^{l}%
\]
suitably reduced modulo $p.$ The coefficient of $\alpha^{n}$ of $\pi_{k,l}$ is
denoted by $\left[  \alpha^{n}\right]  \pi_{k,l}$ and for a general polynumber
$\pi$ we also write%
\[
\left(  \sum_{i=1}^{l}\left[  \alpha^{k_{i}}\right]  \right)  \pi\equiv\left[
\alpha^{k_{1}}\right]  \pi+\left[  \alpha^{k_{2}}\right]  \pi+\cdots+\left[
\alpha^{k_{l}}\right]  \pi.
\]

The key point that connects these circular polynumbers and the super Catalan
numbers is the following restatement of the main result of Georgiadis,
Munemasa, and Tanaka in \cite{Georgiadis} in terms of $\Omega\left(
m,n\right)  $ and $\pi_{k,l}$.

\begin{theorem}
\label{Thm algebraic derivation Munemasa identity intro}For any natural
numbers $m$ and $n$,%
\[
\Omega\left(  m,n\right)  =\left(  -1\right)  ^{n}\left[  \alpha^{m+n}\right]
\pi_{2m,2n}.
\]

\end{theorem}

The complete summation formulas for $\psi_{c,q}$ are presented separately as
Theorems \ref{Thm CIF green group structure}, \ref{Theorem formula red CIF},
and \ref{Thm formula blue CIF combined}, which we summarize below. They
involve geometric constants that depend on the sizes of the unit circles in
each geometry. As usual, $\left\lfloor x\right\rfloor $ denotes the usual
floor function.

\begin{theorem}
\label{Thm CIF formula full Intro}Consider the finite field $\mathbb{F}_{q}$
of size $q=p^{r}$ with $p>2$.

\begin{enumerate}
\item In the green geometry,
\[
\psi_{g,q}\left(  \alpha^{k}\beta^{l}\right)  =\left\{
\begin{array}
[c]{cl}%
1_{q} & \text{if }\left\vert S_{g,q}\right\vert ~|~\left(  k-l\right)
\text{,}\\
0 & \text{otherwise.}%
\end{array}
\right.
\]

\item In the red geometry, with $w\equiv\frac{1}{2}\left\vert S_{r,q}%
\right\vert $ and $R\equiv\left\lfloor \frac{k+l}{\left\vert S_{r,q}%
\right\vert }\right\rfloor $,
\[
\psi_{r,q}(\alpha^{k}\beta^{l})=\left\{
\begin{array}
[c]{cl}%
\left(
{\displaystyle\sum\limits_{\left\vert d\right\vert \leq R}}
\left[  \alpha^{m+n+dw}\right]  \right)  \pi_{k,l}\operatorname{mod}p &
\text{if }k=2m\text{, }l=2n\text{,}\\
0 & \text{otherwise.}%
\end{array}
\right.
\]

\item In the blue geometry, with $w\equiv\frac{1}{2}\left\vert S_{b,q}%
\right\vert $ and $R\equiv\left\lfloor \frac{k+l}{\left\vert S_{b,q}%
\right\vert }\right\rfloor $,
\[
\psi_{b,q}(\alpha^{k}\beta^{l})=\left\{
\begin{array}
[c]{cl}%
\left(
{\displaystyle\sum\limits_{\left\vert d\right\vert \leq R}}
\left[  \alpha^{m+n+dw}\right]  \right)  \left(  -1\right)  ^{n}\pi
_{k,l}\operatorname{mod}p & \text{if }k=2m\text{, }l=2n\text{,}\\
0 & \text{otherwise.}%
\end{array}
\right.
\]

\end{enumerate}
\end{theorem}

Thus we see that the complete formulas to evaluate $\psi_{r,q}(\alpha^{k}%
\beta^{l})$ and $\psi_{b,q}(\alpha^{k}\beta^{l})$ involve an equally-spaced
ladder of coefficients of $\pm\pi_{k,l}$ symmetrically placed around the
central coefficient which is, by the result of Georgiadis, Munemasa, and
Tanaka, $\Omega\left(  m,n\right)  $ when $k=2m$ and $l=2n$.

To illustrate, we will compute the value of $\psi_{b,q}$ evaluated at
$\alpha^{6}\beta^{2}$ over two different fields $\mathbb{F}_{13}$ and
$\mathbb{F}_{7}$. Both computations involve certain coefficients of the
polynumber
\begin{align*}
-\pi_{6,2} &  =-\left(  \frac{1+\alpha}{2}\right)  ^{6}\left(  \frac{1-\alpha
}{2}\right)  ^{2}\\
&  =-\frac{1}{256}\alpha^{8}-\frac{1}{64}\alpha^{7}-\frac{1}{64}\alpha
^{6}+\frac{1}{64}\alpha^{5}+\frac{5}{128}\alpha^{4}+\frac{1}{64}%
\allowbreak\alpha^{3}-\frac{1}{64}\alpha^{2}-\frac{1}{64}\alpha-\frac{1}{256}.
\end{align*}

Note that this expression is symmetric around the central term of $\alpha
^{4},$ where the coefficient is $\frac{5}{128}=\Omega\left(  3,1\right)  .$

\begin{example}
\label{Example q = 13 continued}Continuing our previous example where $q=13$,
we deduce that $w=\frac{1}{2}\left\vert S_{b,13}\right\vert =6$ and
$R_{k,l}=\left\lfloor \frac{k+l}{\left\vert S_{b,13}\right\vert }\right\rfloor
=0.$ The formula for $\psi_{b,q}$ in Theorem \ref{Thm CIF formula full Intro}
above then reads
\[
\psi_{b,13}(\alpha^{6}\beta^{2})=\left(
{\displaystyle\sum\limits_{\left\vert d\right\vert \leq0}}
\left[  \alpha^{4+6d}\right]  \right)  \left(  -\pi_{8,2}\right)
\operatorname{mod}13=\left[  \alpha^{4}\right]  \left(  -\pi_{8,2}\right)
\operatorname{mod}13=4.
\]
Here we see that the condition $0\leq k+l<q-1$ reduces the sum to just the
central term.

However when $q=7$, $-1$ is not a square so by Lemma
\ref{Lemma size of unit circles Intro}, $w=\frac{1}{2}\left\vert
S_{b,7}\right\vert =4$ and $R_{k,l}=\left\lfloor \frac{k+l}{\left\vert
S_{b,7}\right\vert }\right\rfloor =1.$ It follows that
\begin{align*}
\psi_{b,7}(\alpha^{6}\beta^{2}) &  =\left(
{\displaystyle\sum\limits_{\left\vert d\right\vert \leq1}}
\left[  \alpha^{4+4d}\right]  \right)  \left(  -\pi_{8,2}\right)
\operatorname{mod}7\\
&  =\left(  \left[  \alpha^{0}\right]  +\left[  \alpha^{4}\right]  +\left[
\alpha^{8}\right]  \right)  \left(  -\pi_{8,2}\right)  \operatorname{mod}7\\
&  =\left(  -\frac{1}{256}+\frac{5}{128}-\frac{1}{256}\right)
\operatorname{mod}7=\frac{1}{32}\operatorname{mod}7=2.
\end{align*}
It can be verified that $S_{b,7}=\left\{  \left[  1,0\right]  \ \left[
6,0\right]  \ \left[  0,1\right]  \ \left[  0,1\right]  \ \left[  2,2\right]
\ \left[  5,5\right]  \ \left[  2,5\right]  \ \left[  5,2\right]  \right\}  $
and by manual computation $\psi_{b,7}(\alpha^{6}\beta^{2})=2^{6}2^{2}%
+5^{6}5^{2}+2^{6}5^{2}+5^{6}2^{2}=2.$
\end{example}

\subsection{Connection to characteristic zero integration theory}

A formula for the integral of a polynomial over a general sphere was given by
Baker \cite{Baker} and later by Folland \cite{Folland}, both of which involve
the gamma function $\Gamma$. Here we restate the main theorem in Folland's
paper involving the usual rotationally-invariant measure $d\mu$ on the
$\left(  d-1\right)  $-dimensional unit sphere in $\mathbb{R}^{d}$ coming from
the standard Riemannian structure on $\mathbb{R}^{d}$, without proof.

\begin{theorem}
\label{Thm Folland formula Intro}Let $S_{d}=\left\{  \left[  x_{1}%
,x_{2},\ldots,x_{d}\right]  \in\mathbb{R}^{d}\ \colon\ x_{1}^{2}+x_{2}%
^{2}+\cdots+x_{d}^{2}=1\right\}  $. We have that
\[
\int_{S_{d}}x_{1}^{a_{1}}x_{2}^{a_{2}}\cdots x_{d}^{a_{d}}~d\mu=\left\{
\begin{array}
[c]{cl}%
\dfrac{2}{\Gamma\left(  b_{1}+b_{2}+\cdots+b_{d}\right)  }%
{\displaystyle\prod\limits_{i=1}^{d}}
\Gamma\left(  b_{i}\right)  & \text{if }\left(  a_{1},a_{2},\ldots
,a_{n}\right)  =\left(  2b_{1},2b_{2},\ldots,2b_{n}\right) \\
0 & \text{otherwise.}%
\end{array}
\right.
\]

\end{theorem}

Importantly, we note that the formula in Theorem
\ref{Thm Folland formula Intro} above can be reformulated without reference to
$\Gamma$ and $\pi$, thereby eliminating the transcendental aspect, but the key
is to renormalize the integral. We illustrate this in the specific case $d=2$
involving the usual (blue) unit circle. Indeed, by noting that $\Gamma\left(
x+1\right)  =x\Gamma\left(  x\right)  $ and $\Gamma\left(  \frac{1}{2}\right)
=\allowbreak\sqrt{\pi}$, we obtain a formula for the integral when $\left(
a_{1},a_{2}\right)  =\left(  2m,2n\right)  $ as
\begin{align*}
\frac{2\Gamma\left(  b_{1}\right)  \Gamma\left(  b_{2}\right)  }{\Gamma\left(
b_{1}+b_{2}\right)  }  &  =\frac{2\Gamma\left(  m+\frac{1}{2}\right)
\Gamma\left(  n+\frac{1}{2}\right)  }{\Gamma\left(  m+n+1\right)  }\\
&  =\frac{2}{\left(  m+n\right)  !}\left(
{\displaystyle\prod\limits_{k=1}^{m}}
\left(  \frac{2k-1}{2}\right)  \right)  \left(
{\displaystyle\prod\limits_{l=1}^{n}}
\left(  \frac{2l-1}{2}\right)  \right)  \Gamma\left(  \frac{1}{2}\right)
^{2}\\
&  =\frac{2\pi}{2^{m+n}\left(  m+n\right)  !}\frac{\left(  2m\right)  !}%
{2^{m}m!}\frac{\left(  2n\right)  !}{2^{n}n!}\\
&  =2\pi\Omega\left(  m,n\right)  .
\end{align*}

If the integral over the unit circle $S_{2}$ is normalized, then
\begin{equation}
\int_{S_{2}}x_{1}^{a_{1}}x_{2}^{a_{2}}~d\mu=\left\{
\begin{array}
[c]{cl}%
\Omega(m,n) & \text{if }a_{1}=2m\text{, }a_{2}=2n\text{,}\\
0 & \text{otherwise}%
\end{array}
\right.  \nonumber%
\end{equation}
which is strikingly similar to the formula for $\psi_{b,q}$ presented in
Theorem \ref{Thm CIF formula SC Intro}. 

\section{Polynumbers and functions over finite fields}

Since we are dealing here with finite mathematics and finite arithmetic, we
propose to dispense with the usual reliance on "infinite sets". Throughout
this paper, $\mathbb{N}$ will denote the \textit{type} of a natural number ---
which for us includes the number $0$; $\mathbb{Z}$ denotes the type of
integers; $\mathbb{Q}$ denotes the type of rational numbers; and
$\mathbb{F}_{q}$ will denote the finite field of size $q,$ where $q$ is a
power of an odd prime $p$. We may express the multiplicative identity in
$\mathbb{F}_{q}$ as $1_{q}$ if there is possibility of confusion. The notation
$\mathbb{F}^{\times}$ is reserved for the multiplicative group of the field
$\mathbb{F}$. We also use the convention that $0^{0}=1$.

We briefly distinguish between polynomials as formal expressions and as
functions by regarding the former purely as finite lists of numbers in
$\mathbb{F}$, thereby avoiding the notion of "variables". So a
\textbf{polynumber} $\pi$ over a field $\mathbb{F}$ is a finite array of
numbers in $\mathbb{F}$ written vertically with a bar above it and another one
to the left, as in
\[
\pi=\left\vert \overline{%
\begin{array}
[c]{c}%
u_{0}\\
u_{1}\\
\vdots\\
u_{k}%
\end{array}
}\right.  .
\]
Adding or removing zeroes from the bottom will not change a polynumber. The
\textbf{zero polynumber} and the\textbf{ identity polynumber} are respectively
$\mathbf{0}\equiv\left\vert \overline{%
\begin{array}
[c]{c}%
0
\end{array}
}\right.  $ and $\mathbf{1}\equiv\left\vert \overline{%
\begin{array}
[c]{c}%
1
\end{array}
}\right.  $. The addition and scalar multiplication of polynumbers are
entry-wise while polynumber multiplication is defined as the usual Cauchy
product on sequence of numbers. The notion of degree and irreducibility of a
polynumber is similar to that of a polynomial, with the degree of $\mathbf{0}$
not defined. We denote the algebra of polynumbers over $\mathbb{F}$ by
$\mathrm{Pol}\left(  \mathbb{F}\right)  $.

After we define the special polynumber
\[
\alpha\equiv\left\vert \overline{%
\begin{array}
[c]{c}%
0\\
1
\end{array}
}\right.
\]
we say a \textbf{polynomial} in $\alpha$ is a polynumber written as a linear
combination of powers of the polynumber $\alpha$. If $\pi=u_{0}\mathbf{1}%
+u_{1}\alpha+\cdots+u_{k}\alpha^{k}$ is a polynomial in $\alpha$, then $u_{k}$
is the coefficient of $\alpha^{k}$ in $\pi,$ and we write $u_{k}=$ $\left[
\alpha^{k}\right]  \pi$. We will allow ourselves to use also summation
notation and write
\[
\left(  \sum_{i=1}^{l}\left[  \alpha^{k_{i}}\right]  \right)  \pi\equiv\left[
\alpha^{k_{1}}\right]  \pi+\cdots+\left[  \alpha^{k_{l}}\right]  \pi.
\]

In a similar fashion, we can define a $2$-\textbf{polynumber} over
$\mathbb{F}$ to be a two-dimensional array of numbers in $\mathbb{F}$ of the
form%
\[
\pi\equiv\left\vert \overline{%
\begin{array}
[c]{cccc}%
u_{00} & u_{01} & \cdots & u_{0n}\\
u_{10} & u_{11} & \cdots & u_{1n}\\
\vdots & \vdots & \ddots & \vdots\\
u_{m0} & u_{m1} & \cdots & u_{mn}%
\end{array}
}\right.
\]
with the convention that adding further zeroes to the right and down leaves
the $2$-polynumber unchanged. The algebra operations are analogous to those in
$\mathrm{Pol}\left(  \mathbb{F}\right)  $, with the multiplication defined by
the two-dimensional version of the Cauchy product. The algebra of
$2$-polynumbers over $\mathbb{F}$ is denoted by $\mathrm{Pol}_{2}\left(
\mathbb{F}\right)  $. The multiplicative identity is also denoted by
$\mathbf{1}$, not to be confused with that of $\mathrm{Pol}\left(
\mathbb{F}\right)  $ if the context is clear. After we define the polynumber%
\[
\beta\equiv\left\vert \overline{%
\begin{array}
[c]{cc}%
0 & 1
\end{array}
}\right.
\]
any $2$-polynumber can be written as a linear combination of powers of
$\alpha$ and $\beta$ and then the algebra operations agree with the familiar
operations on polynomials. Note however that here and throughout $\alpha$ and
$\beta$ are \textit{not} variables; they are just particular polynumbers.
Clearly we may regard $\mathrm{Pol}\left(  \mathbb{F}\right)  $ as a
subalgebra of $\mathrm{Pol}_{2}\left(  \mathbb{F}\right)  $.

If $\pi\in\mathrm{Pol}\left(  \mathbb{F}_{p}\right)  $ is an irreducible
polynumber of degree $r$ and $\left(  \pi\right)  $ denotes the ideal of
$\mathrm{Pol}\left(  \mathbb{F}_{p}\right)  $ generated by $\pi$, then
$\mathrm{Pol}\left(  \mathbb{F}_{p}\right)  /\left(  \pi\right)  $ is a field
of size $q=p^{r}$. Since any two finite fields of the same size are
isomorphic, the choice of $\pi$ does not matter.

The two-dimensional \textbf{affine plane} is denoted by $\mathbb{A=A}\left(
\mathbb{F}\right)  =\left\{  \left[  x,y\right]  \ \colon\ x,y\in
\mathbb{F}\right\}  $, with the objects $\left[  a,b\right]  $ called
\textbf{affine points}, or just \textbf{points}. Moreover, the space of
functions from $\mathbb{A}$ to $\mathbb{F}$, denoted by $\mathbb{F}%
^{\mathbb{A}}$, is an algebra over $\mathbb{F}$ under the usual pointwise
operation. For brevity, we write $f\left(  \left[  x,y\right]  \right)  $
simply as $f\left(  x,y\right)  $.

If $\mathbb{F=F}_{q}$ is a finite field with $q=p^{r}$ elements, then
$\dim\left(  \mathbb{F}_{q}^{\mathbb{A}}\right)  =q^{2}$ and a basis for
$\mathbb{F}_{q}^{\mathbb{A}}$ is $\Delta=\left\{  \delta_{\left[  a,b\right]
}\ \colon\ \left[  a,b\right]  \in\mathbb{A}\right\}  $, where
\[
\delta_{\left[  a,b\right]  }\left(  x,y\right)  =\left\{
\begin{array}
[c]{cl}%
1_{q} & \text{if }\left[  x,y\right]  =\left[  a,b\right] \\
0 & \text{if }\left[  x,y\right]  \neq\left[  a,b\right]  .
\end{array}
\right.
\]
The multiplicative identity in $\mathbb{F}_{q}^{\mathbb{A}}$ is $1_{\mathbb{F}%
^{\mathbb{A}}}$. Since $\delta_{\left[  a,b\right]  }\left(  x,y\right)  $
vanishes at every point in $\mathbb{A}$ except $\left[  a,b\right]  $, we have%
\begin{equation}
1_{\mathbb{F}^{\mathbb{A}}}=\sum_{\left[  a,b\right]  \in\mathbb{A}}%
\delta_{\left[  a,b\right]  }. \label{Property of constant 1 function}%
\end{equation}

Evaluation of a $2$-polynumber in $\mathrm{Pol}_{2}\left(  \mathbb{F}%
_{q}\right)  $ then gives a linear map $\varepsilon$~$\colon\ \mathrm{Pol}%
_{2}\left(  \mathbb{F}_{q}\right)  \rightarrow\mathbb{F}_{q}^{\mathbb{A}}$.
Although clearly not injective, this map is onto. To see this, we set for any
$\left[  a,b\right]  \in\mathbb{A}$ the $2$-polynumber%
\begin{equation}
\pi_{\left[  a,b\right]  }\equiv\left(  \prod_{\substack{t\in\mathbb{F}_{q}\\t\neq
a}}\left(  t-a\right)  \prod_{\substack{u\in\mathbb{F}_{q}\\u\neq b}}\left(
u-b\right)  \right)  ^{-1}\prod_{\substack{t\in\mathbb{F}_{q}\\t\neq
a}}\left(  t-\alpha\right)  \prod_{\substack{u\in\mathbb{F}_{q}\\u\neq
b}}\left(  u-\beta\right)  \nonumber%
\end{equation}
and observe that $\varepsilon\left(  \pi_{\left[  a,b\right]  }\right)
=\delta_{\left[  a,b\right]  }$.

The subspace $\mathcal{P}_{q}=\operatorname*{span}\left\{  \alpha^{k}\beta
^{l}\ \colon\ 0\leq k,l\leq q-1\right\}  $ of $\mathrm{Pol}_{2}\left(
\mathbb{F}_{q}\right)  $ has dimension $q^{2}$ and we call it the
\textbf{principal subspace}. Note that for any point $\left[  a,b\right]  $ we
have $\pi_{\left[  a,b\right]  }\in\mathcal{P}_{q}$. Since $\dim\left(
\mathcal{P}_{q}\right)  =\dim\left(  \mathbb{F}_{q}^{\mathbb{A}}\right)  $ and
$\left.  \varepsilon\right\vert _{\mathcal{P}_{q}}$ is onto, it is a vector
space isomorphism. In particular, the set $\Pi=\left\{  \pi_{\left[
a,b\right]  }\ \colon\ \left[  a,b\right]  \in\mathbb{A}\right\}  $ is a basis
for $\mathcal{P}_{q}$. The following result gives a vector space decomposition
of $\mathrm{Pol}_{2}\left(  \mathbb{F}_{q}\right)  $.

\begin{proposition}
\label{Decomposition of Pol space}We have that $\mathrm{Pol}_{2}\left(
\mathbb{F}_{q}\right)  =\mathcal{P}_{q}\oplus\ker\left(  \varepsilon\right)  $.
\end{proposition}

\begin{proof}
Take any $\pi\in\mathrm{Pol}_{2}\left(  \mathbb{F}_{q}\right)  $. Since
$\left.  \varepsilon\right\vert _{\mathcal{P}_{q}}$ is an isomorphism, there
is a unique $\pi_{1}=\left(  \left.  \varepsilon\right\vert _{\mathcal{P}_{q}%
}^{-1}\circ\varepsilon\right)  \left(  \pi\right)  \in\mathcal{P}_{q}$ for
which $\varepsilon\left(  \pi_{1}\right)  =\varepsilon\left(  \pi\right)  $.
By linearity, $\pi-\pi_{1}\in\ker\left(  \varepsilon\right)  $ and therefore
the decomposition $\pi_{1}+\left(  \pi-\pi_{1}\right)  $ is unique.
\end{proof}

\subsection{Group actions on various spaces}

We now define some group actions on the relevant spaces that will be useful
later in our analysis. The traditional treatment to specifying a group action
is by fixing a direction the group acts from, that is, agreeing that every
action will be a left action or a right action. We will instead be more
flexible and specify explicitly whether a given action is a left action or a
right action.

We say that a group $G$ with identity $e$ \textbf{left acts} on a space
$\mathbb{X}$ precisely when for any $g\in G$ and $X\in\mathbb{X}$ there is an
element $g\cdot X\in\mathbb{X}$ such that $e\cdot X=X$ and $h\cdot\left(
g\cdot X\right)  =\left(  hg\right)  \cdot X$ for $X\in\mathbb{X}$ and $g,h\in
G.$ We say that a group $G$ with identity $e$ \textbf{right acts} on a space
$\mathbb{X}$ precisely when for any $g\in G$ and $X\in\mathbb{X}$ there is an
element $X\cdot g\in\mathbb{X}$ such that $X\cdot e=X$ and $\left(  X\cdot
g\right)  \cdot h=X\cdot\left(  gh\right)  $ for $X\in\mathbb{X}$ and $g,h\in
G$.

From our point of view, the distinction between these two types of actions has
at least two benefits. First, notions of duality and functiorality can
generally be more easily appreciated. Secondly, this approach generally avoids
the need for inverses, which is convenient since they are generally
computationally expensive and thus should be avoided theoretically whenever possible.

Now suppose that $G$ is a group with identity $e$ that right acts on the
affine plane $\mathbb{A}$. Then there is a natural left action of $G$ on the
space of functions $\mathbb{F}_{q}^{\mathbb{A}}$ defined for $h\in G$ and
$f\in\mathbb{F}_{q}^{\mathbb{A}}$ by $\left(  h\cdot f\right)  \left(
x,y\right)  \equiv f\left(  \left[  x,y\right]  \cdot h\right)  $. We can
quickly verify that the basis $\Delta$ is invariant under this left action:
$\left(  h\cdot\delta_{\left[  a,b\right]  }\right)  \left(  x,y\right)
=\delta_{\left[  a,b\right]  }\left(  \left[  x,y\right]  \cdot h\right)  $
for any $h\in G$ and $\left[  a,b\right]  \in\mathbb{A}$, and therefore
\begin{equation}
h\cdot\delta_{\left[  a,b\right]  }=\delta_{\left[  a,b\right]  \cdot h^{-1}}.
\label{Action of G to delta function}%
\end{equation}

We say that $f\in\mathbb{F}_{q}^{\mathbb{A}}$ is $G$-\textbf{invariant}
precisely when $h\cdot f=f$ ~for any $h\in G$, and let $\mathcal{I}_{G}$ be
the subalgebra of $\mathbb{F}_{q}^{\mathbb{A}}$ consisting of $G$-invariant
functions. Note that for any $f\in\mathbb{F}_{q}^{\mathbb{A}}$, the
$G$-\textbf{average} of $f$, defined by
\[
\rho_{G}\left(  f\right)  \equiv\frac{1}{\left\vert G\right\vert }\sum_{h\in
G}\left(  h\cdot f\right)
\]
is always $G$-invariant, provided that $\left\vert G\right\vert \neq0$ in
$\mathbb{F}_{q}$. Indeed, for any $g\in G$ and $\left[  x,y\right]
\in\mathbb{A}$, we have
\begin{align*}
\left(  g\cdot\rho_{G}\left(  f\right)  \right)  \left(  x,y\right)   &
=\rho_{G}\left(  f\right)  \left(  \left[  x,y\right]  \cdot g\right) \\
&  =\frac{1}{\left\vert G\right\vert }\sum_{h\in G}\left(  h\cdot f\right)
\left(  \left[  x,y\right]  \cdot g\right)  =\frac{1}{\left\vert G\right\vert
}\sum_{h\in G}\left(  gh\cdot f\right)  \left(  x,y\right) \\
&  =\frac{1}{\left\vert G\right\vert }\sum_{h\in G}\left(  h\cdot f\right)
\left(  x,y\right)  \text{ (via a change of index }h\mapsto g^{-1}h\text{)}\\
&  =\rho_{G}\left(  f\right)  \left(  x,y\right)  .
\end{align*}
Moreover, it is obvious that $\rho_{G}\left(  f\right)  =f$ if $f\in
\mathcal{I}_{G}$. This makes $\rho_{G}$ a projection onto $\mathcal{I}_{G}.$

Now there is also a natural right action of $G$ on the dual space $\left(
\mathbb{F}_{q}^{\mathbb{A}}\right)  ^{\ast}$ of $\mathbb{F}_{q}^{\mathbb{A}}$
by $\left(  \phi\cdot h\right)  \left(  f\right)  \equiv\phi\left(  h\cdot
f\right)  $ for any $\phi\in\left(  \mathbb{F}_{q}^{\mathbb{A}}\right)
^{\ast}$, $h\in G$, and $f\in\mathbb{F}_{q}^{\mathbb{A}}$. Here the duality
between $\mathbb{F}_{q}^{\mathbb{A}}$ and $\left(  \mathbb{F}_{q}^{\mathbb{A}%
}\right)  ^{\ast}$ is also reflected in the left and right action of $G$ on
$\mathbb{F}_{q}^{\mathbb{A}}$ and $\left(  \mathbb{F}_{q}^{\mathbb{A}}\right)
^{\ast}$.

A linear functional $\phi\in\left(  \mathbb{F}_{q}^{\mathbb{A}}\right)
^{\ast}$ is $G$-\textbf{invariant }precisely when $\phi\cdot h=\phi$ for any
$h\in G$. The subspace of $G$-invariant functionals in $\left(  \mathbb{F}%
_{q}^{\mathbb{A}}\right)  ^{\ast}$ can then be identified with $\mathcal{I}%
_{G}^{\ast}$. By a similar argument, the $G$-\textbf{average} of $\phi$
defined by%
\[
\widehat{\rho}_{G}\left(  \phi\right)  \equiv\frac{1}{\left\vert G\right\vert
}\sum_{h\in G}\left(  \phi\cdot h\right)
\]
provided that $\left\vert G\right\vert \neq0$ in $\mathbb{F}_{q}$ is
$G$-\textbf{invariant} and therefore $\widehat{\rho}_{G}$ is a projection of
$\left(  \mathbb{F}_{q}^{\mathbb{A}}\right)  ^{\ast}$ onto $\mathcal{I}%
_{G}^{\ast}$.

\subsection{Krawtchouk polynumbers and the circular polynumbers}

For $l\in\mathbb{N}$, we define the rational polynumber
\[
\dbinom{\alpha}{l}\equiv\left\{
\begin{array}
[c]{cl}%
\dfrac{\alpha\left(  \alpha-1\right)  \left(  \alpha-2\right)  \ldots\left(
\alpha-l+1\right)  }{l!} & \text{if }l\geq1\\
\mathbf{1} & \text{if }l=0.
\end{array}
\right.
\]
With this notation, the composition of $\dbinom{\alpha}{l}$ for $l\geq1$ with
another rational polynumber $\pi$ can be written as%
\[
\dbinom{\pi}{l}=\frac{\pi\left(  \pi-1\right)  \left(  \pi-2\right)
\ldots\left(  \pi-l+1\right)  }{l!}.
\]
For $0\leq n\leq d$, the $n$-th \textbf{Krawtchouk polynumber} $k_{n}^{\left(
d\right)  }$ of order $d$ is given by%
\[
k_{n}^{\left(d\right)}\left(\alpha\right)  \equiv\sum_{l=0}^{n}\left(
-1\right)  ^{l}\dbinom{\alpha}{l}\dbinom{d-\alpha}{n-l}.
\]
This family of polynumbers generalizes the Hermite polynumbers and was
introduced by Mikhail Krawtchouk in the late 1920s \cite{Krawtchouk}. Their
generating function is particularly pleasant, given by%
\begin{equation}
\sum_{n=0}^{d}k_{n}^{\left(d\right)}\left(m\right)  \alpha^{n}=\left(
1-\alpha\right)  ^{m}\left(  1+\alpha\right)  ^{d-m}
\label{GF Krawtchouk polynumbers}%
\end{equation}
where $k_{n}^{\left(  d\right)  }\left(  m\right)  $ is the evaluation of
$k_{n}^{\left(  d\right)  }$ at $m\in\left\{  0,1,2,\ldots,d\right\}  $.

Krawtchouk polynumbers play an important role in many areas of mathematics,
such as harmonic analysis in the context of orthogonal polynomials
\cite{Dunkl}, \cite{Dunkl2}, coding theory in the context of bounds of codes
\cite{Levensthein}, \cite{Ling}, \cite{MacWilliams}, and probability theory in
the context of symmetric random walk \cite{Feinsilver2}.

In 2011, Georgiadis, Munemasa, and Tanaka in their paper \cite{Georgiadis}
gave an interpretation of $S\left(  m,n\right)  $ in terms of values of
Krawtchouk polynumbers. We state their main result now, and will provide a new
proof of it shortly.

\begin{theorem}
[Georgiadis, Munemasa, Tanaka]\label{Thm Georgiadis etc}For any natural
numbers $m$ and $n$,
\begin{equation}
k_{m+n}^{\left(  2m+2n\right)  }\left(  2m\right)  =\left(  -1\right)
^{m}S\left(  m,n\right)  . \label{Georgiadis etc original result}%
\end{equation}

\end{theorem}

The next proposition is a restatement of Theorem \ref{Thm Georgiadis etc} in
terms of the generating function.

\begin{proposition}
\label{Prop SC numbers and GF Krawtchouk}For any natural numbers $m$ and $n$,%
\[
S\left(  m,n\right)  =\left(  -1\right)  ^{n}\left[  \alpha^{m+n}\right]
\left(  1+\alpha\right)  ^{2m}\left(  1-\alpha\right)  ^{2n}.
\]

\end{proposition}

\begin{proof}
From (\ref{GF Krawtchouk polynumbers}) and
(\ref{Georgiadis etc original result}), we deduce that
\begin{equation}
\left(  -1\right)  ^{m}S\left(  m,n\right)  =\left[  \alpha^{m+n}\right]
\left(  1-\alpha\right)  ^{2m}\left(  1+\alpha\right)  ^{2n}\text{.}
\label{Cor GF Georgiadis}%
\end{equation}
Finally, by multiplying both sides of equation (\ref{Cor GF Georgiadis}) above by
$\left(  -1\right)  ^{m}$, we obtain%
\begin{align*}
S\left(  m,n\right)   &  =\left(  -1\right)  ^{m}\left[  \alpha^{m+n}\right]
\left(  1-\alpha\right)  ^{2m}\left(  1+\alpha\right)  ^{2n}\\
&  =\left(  -1\right)  ^{n}\left[  \left(  -\alpha\right)  ^{m+n}\right]
\left(  1-\alpha\right)  ^{2m}\left(  1+\alpha\right)  ^{2n}\\
&  =\left[  \alpha^{m+n}\right]  \left(  1+\alpha\right)  ^{2m}\left(
1-\alpha\right)  ^{2n}. \qedhere
\end{align*}
\end{proof}

We now introduce a central family of polynumbers called the \textbf{circular
polynumbers }%
\begin{equation}
\pi_{k,l}\equiv\left(  \frac{1+\alpha}{2}\right)  ^{k}\left(  \frac{1-\alpha
}{2}\right)  ^{l} \nonumber%
\end{equation}
for $k,l\in\mathbb{N}$. These polynumbers, as we shall soon see, are crucial
in describing the Fourier summation functional formulas in the red and blue geometries.

It is clear that $\pi_{0,0}=\mathbf{1}$. Here is a table of the polynumbers
$\pi_{k,l}$ for the next few values of $k$ and $l$, organized by the sum
$k+l$:%
\[
{\normalsize
\begin{tabular}
[c]{c|c}%
$\pi_{k,l}$ & $k+l=1$\\\hline
$\pi_{0,1}$ & ${\scriptsize \allowbreak}\frac{1}{2}\left(  1-\alpha\right)
$\\
$\pi_{1,0}$ & ${\scriptsize \allowbreak}\frac{1}{2}\left(  1+\alpha\right)  $%
\end{tabular}
\ \ }\quad%
\begin{tabular}
[c]{c|c}%
${\normalsize \pi}_{k,l}$ & ${\normalsize k+l=2}$\\\hline
${\normalsize \pi}_{0,2}$ & $\frac{1}{4}\left(  1-2\alpha+\alpha^{2}\right)
$\\
${\normalsize \pi}_{1,1}$ & $\frac{1}{4}\left(  1-\alpha^{2}\right)  $\\
${\normalsize \pi}_{2,0}$ & $\frac{1}{4}\left(  1+2\alpha+\alpha^{2}\right)  $%
\end{tabular}
\ \ \quad%
\begin{tabular}
[c]{c|c}%
${\normalsize \pi}_{k,l}$ & ${\normalsize k+l=3}$\\\hline
${\normalsize \pi}_{0,3}$ & $\frac{1}{8}\left(  1-3\alpha+3\alpha^{2}%
-\alpha^{3}\right)  $\\
$\pi_{1,2}$ & $\frac{1}{8}\left(  1-\alpha-\alpha^{2}+\alpha^{3}\right)  $\\
${\normalsize \pi}_{2,1}$ & $\frac{1}{8}\left(  1+\alpha-\alpha^{2}-\alpha
^{3}\right)  $\\
${\normalsize \pi}_{3,0}$ & $\frac{1}{8}\left(  1+3\alpha+3\alpha^{2}%
+\alpha^{3}\right)  $%
\end{tabular}
\ \ \ \ \ \ \ \ \ \
\]%
\[%
\begin{tabular}
[c]{c|c}%
${\normalsize \pi}_{k,l}$ & ${\normalsize k+l=4}$\\\hline
${\normalsize \pi}_{0,4}$ & $\frac{1}{16}\left(  1-4\alpha+6\alpha^{2}%
-4\alpha^{3}+\alpha^{4}\right)  $\\
${\normalsize \pi}_{1,3}$ & $\frac{1}{16}\left(  1-2\alpha+2\alpha^{3}%
-\alpha^{4}\right)  $\\
${\normalsize \pi}_{2,2}$ & $\frac{1}{16}\left(  1-2\alpha^{2}+\alpha
^{4}\right)  $\\
${\normalsize \pi}_{3,1}$ & $\frac{1}{16}\left(  1+2\alpha-2\alpha^{3}%
-\alpha^{4}\right)  $\\
$\pi_{4,0}$ & $\frac{1}{16}\left(  1+4\alpha+6\alpha^{2}+4\alpha^{3}%
+\alpha^{4}\right)  $%
\end{tabular}
\ \ \quad%
\begin{tabular}
[c]{c|c}%
${\normalsize \pi}_{k,l}$ & ${\normalsize k+l=5}$\\\hline
${\normalsize \pi}_{0,5}$ & $\frac{1}{32}\left(  1-5\alpha+10\alpha
^{2}-10\alpha^{3}+5\alpha^{4}-\alpha^{5}\right)  $\\
${\normalsize \pi}_{1,4}$ & $\frac{1}{32}\left(  1-3\alpha+2\alpha^{2}%
+2\alpha^{3}-3\alpha^{4}+\alpha^{5}\right)  $\\
${\normalsize \pi}_{2,3}$ & $\frac{1}{32}\left(  1-\alpha-2\alpha^{2}%
+2\alpha^{3}+\alpha^{4}-\alpha^{5}\right)  $\\
${\normalsize \pi}_{3,2}$ & $\frac{1}{32}\left(  1+\alpha-2\alpha^{2}%
-2\alpha^{3}+\alpha^{4}+\alpha^{5}\right)  $\\
$\pi_{4,1}$ & $\frac{1}{32}\left(  1+3\alpha+2\alpha^{2}-2\alpha^{3}%
-3\alpha^{3}-\alpha^{5}\right)  $\\
$\pi_{5,0}$ & $\frac{1}{32}\left(  1+5\alpha+10\alpha^{2}+10\alpha^{3}%
+5\alpha^{4}+\alpha^{5}\right)  $%
\end{tabular}
\ \ \ \ \ \ \ \ \ \
\]

\begin{center}
\textbf{Table 3}: The circular polynumbers $\pi_{k,l}$ for $0\leq k+l\leq5$.
\end{center}

We say that a polynumber $\pi=u_{0}\alpha^{0}+u_{1}\alpha^{1}+\cdots
+u_{d-1}\alpha^{d-1}+u_{d}\alpha^{d}$ of degree $d$ is \textbf{palindromic}
precisely when $u_{k}=u_{d-k}$ and \textbf{anti-palindromic} precisely when
$u_{k}=-u_{d-k}$, for all $0\leq k\leq d$. The product of two palindromic
polynumbers or two anti-palindromic polynumbers is palindromic, while the
product of a palindromic and an anti-palindromic polynumber is anti-palindromic.

\begin{lemma}
[Palindromic and symmetry properties of $\pi_{k,l}$]%
\label{Lemma palindromic property Krawtchouk}If $0\leq i,j\leq k+l$ and
$i+j=k+l$ then
\[
\left[  \alpha^{i}\right]  \pi_{k,l}=\left(  -1\right)  ^{l}\left[  \alpha
^{j}\right]  \pi_{k,l}.
\]
Moreover, if $0\leq i\leq k+l$ then
\[
\left[  \alpha^{i}\right]  \pi_{k,l}=\left(  -1\right)  ^{i}\left[  \alpha
^{i}\right]  \pi_{l,k}.
\]

\end{lemma}

\begin{proof}
Clearly $\pi_{k,0}$ is a palindromic polynumber for any $k$, while $\pi_{0,l}$
is palindromic if $l$ is even and anti-palindromic if $l$ is odd. The first
identity follows from $\pi_{k,l}=\pi_{k,0}\pi_{0,l}$. For the second identity,
we observe that $\pi_{l,k}\left(  -\alpha\right)  =\pi_{k,l}\left(
\alpha\right)  $. Thus by extracting the coefficient of $\alpha^{i}$ from both
sides,
\[
\left[  \alpha^{i}\right]  \pi_{k,l}=\left[  \alpha^{i}\right]  \pi
_{l,k}\left(  -\alpha\right)  =\left[  \left(  -\alpha\right)  ^{i}\right]
\pi_{l,k}=\left[  \left(  -1\right)  ^{i}\alpha^{i}\right]  \pi_{l,k}=\left(
-1\right)  ^{i}\left[  \alpha^{i}\right]  \pi_{l,k}. \qedhere
\]
\end{proof}

We now restate Theorem \ref{Thm Georgiadis etc} in terms of the circular super
Catalan numbers $\Omega\left(  m,n\right)  $ and the central coefficient of
the circular polynumbers $\pi_{2m,2n},$ and give an independent algebraic
proof. This allows us to naturally identify the remaining coefficients around
the central one as secondary but still important variants of the circular
super Catalan numbers.

\begin{theorem}
\label{Thm algebraic derivation Munemasa identity}For any natural numbers $m$
and $n$,%
\[
\Omega\left(  m,n\right)  =\left(  -1\right)  ^{n}\left[  \alpha^{m+n}\right]
\pi_{2m,2n}.
\]

\end{theorem}

\begin{proof}
We calculate that%
\begin{align*}
\left(  -1\right)  ^{n}\left[  \alpha^{m+n}\right]  \pi_{2m,2n}  &  =\left(
-1\right)  ^{n}\left[  \alpha^{m+n}\right]  \left(  \left(  \frac{1+\alpha}%
{2}\right)  ^{2m}\left(  \frac{1-\alpha}{2}\right)  ^{2n}\right) \\
&  =\frac{\left(  -1\right)  ^{n}}{4^{m+n}}\sum_{\substack{0\leq
s\leq2m\\0\leq t\leq2n\\s+t=m+n}}\left(  -1\right)  ^{t}\dbinom{2m}{s}%
\dbinom{2n}{t}.
\end{align*}
By using the convention that $\dbinom{k}{s}=0$ if $s<0$ or $s>k$, we may
rearrange this sum as
\begin{align}
\left(  -1\right)  ^{n}\left[  \alpha^{m+n}\right]  \pi_{2m,2n}  &
=\frac{\left(  -1\right)  ^{n}}{4^{m+n}}\sum_{s=0}^{m+n}\left(  -1\right)
^{m+n-s}\dbinom{2m}{s}\dbinom{2n}{m+n-s}\nonumber\\
&  =\frac{\left(  -1\right)  ^{m}}{4^{m+n}}\sum_{s=0}^{m+n}\left(  -1\right)
^{s}\frac{\left(  2m\right)  !}{s!\left(  2m-s\right)  !}\frac{\left(
2n\right)  !}{\left(  m+n-s\right)  !\left(  -m+n+s\right)  !}\nonumber\\
&  =\frac{\left(  -1\right)  ^{m}\left(  2m\right)  !\left(  2n\right)
!}{4^{m+n}\left(  m+n\right)  !\left(  m+n\right)  !}\sum_{s=0}^{m+n}\left(
-1\right)  ^{s}\dbinom{m+n}{s}\dbinom{m+n}{2m-s}.
\label{Algebraic derivation Munemasa identity 2}%
\end{align}
The summation on the right-hand side of
(\ref{Algebraic derivation Munemasa identity 2}) can be realized as
\[
\left[  \alpha^{2m}\right]  \left(  1+\alpha\right)  ^{m+n}\left(
1-\alpha\right)  ^{m+n}=\left[  \alpha^{2m}\right]  \left(  1-\alpha
^{2}\right)  ^{m+n}%
\]
so we get
\begin{align*}
\left(  -1\right)  ^{n}\left[  \alpha^{m+n}\right]  \pi_{2m,2n}  &
=\frac{\left(  -1\right)  ^{m}\left(  2m\right)  !\left(  2n\right)
!}{4^{m+n}\left(  m+n\right)  !\left(  m+n\right)  !}\left[  \alpha
^{2m}\right]  \left(  1-\alpha^{2}\right)  ^{m+n}\\
&  =\frac{\left(  2m\right)  !\left(  2n\right)  !}{4^{m+n}\left(  m+n\right)
!\left(  m+n\right)  !}\dbinom{m+n}{m}=\Omega\left(  m,n\right). \qedhere
\end{align*}
\end{proof}

\section{Chromogeometry and Dihedrons}

We use chromogeometry, a three-fold symmetry introduced by Wildberger in \cite{Wildberger},
as a foundation for Fourier summation over finite fields, involving both the
unit circles over which Fourier summation takes place, and the associated
groups of rotation matrices which act on the space of $2$-polynumbers. The
three geometries will be defined through symmetric bilinear forms induced by
symmetric matrices, thereby avoiding notions of distance involving square
roots, which do not generally apply to finite fields.

If $\mathbb{F}$ is a field of characteristic not two, then the algebra of
$2\times2$ matrices over $\mathbb{F}$, which we call the \textbf{dihedron
algebra} over $\mathbb{F}$, is denoted by $\mathbb{D}=\mathbb{D}\left(
\mathbb{F}\right)  $. This is a four-dimensional algebra over $\mathbb{F}$
with a special basis consisting of the four matrices%
\[
1_{\mathbb{D}}\equiv\left(
\begin{array}
[c]{cc}%
1_{\mathbb{F}} & 0\\
0 & 1_{\mathbb{F}}%
\end{array}
\right)  \quad i\equiv\left(
\begin{array}
[c]{cc}%
0 & 1_{\mathbb{F}}\\
-1_{\mathbb{F}} & 0
\end{array}
\right)  \quad j\equiv\left(
\begin{array}
[c]{cc}%
0 & 1_{\mathbb{F}}\\
1_{\mathbb{F}} & 0
\end{array}
\right)  \quad k\equiv\left(
\begin{array}
[c]{cc}%
1_{\mathbb{F}} & 0\\
0 & -1_{\mathbb{F}}%
\end{array}
\right)  .
\]
The dihedron algebra terminology is motivated by the observation that
$\left\{  \pm1_{\mathbb{D}},\pm i,\pm j,\pm k\right\}  $ is the dihedral group
$D_{4},$ and was coined by Wildberger in his online \textit{Linear Algebra}
lectures series. The situation is analogous to and shares
many parallels with the quaternion group $Q$ also of order $8$ and the
associated four-dimensional algebra $\mathbb{H}$ of quaternions. We can think
of $\mathbb{F}$ as being embedded in $\mathbb{D}$ by making the identification
$t\mapsto t1_{\mathbb{D}}$.

These four special basis matrices satisfy
\begin{equation}%
\begin{array}
[c]{ccc}%
i^{2}=-1_{\mathbb{D}} & j^{2}=1_{\mathbb{D}} & k^{2}=1_{\mathbb{D}}\\
jk=-i & ki=j & ij=k\\
kj=i & ik=-j & ji=-k.
\end{array}
\nonumber%
\end{equation}
The difference between $\mathbb{D}$ and $\mathbb{H}$ then amounts to small
alterations in signs for the structure constants when comparing the basis
elements $1_{\mathbb{D}}$, $i$, $j$, and $k$ and the more usual corresponding
elements in $Q$. We equip $\mathbb{D}$ with the structure of a $\ast
$-\textbf{algebra }over $\mathbb{F}$ by defining the involution to be the
\textbf{adjugate} map
\[
\left(
\begin{array}
[c]{cc}%
a & b\\
c & d
\end{array}
\right)  ^{\ast}\equiv\left(
\begin{array}
[c]{cc}%
d & -b\\
-c & a
\end{array}
\right)
\]
and note that for $h_{1},h_{2}\in\mathbb{D}$ we have $\left(  h_{1}%
h_{2}\right)  ^{\ast}=h_{2}^{\ast}h_{1}^{\ast}$. It is also obvious that
$hh^{\ast}=h^{\ast}h=\det\left(  h\right)  1_{\mathbb{D}}$ for all
$h\in\mathbb{D}$.

There are three particularly useful commutative $\ast$-subalgebras of
$\mathbb{D}\left(  \mathbb{F}\right)  $, namely
\[
\mathbb{C}_{b}\left(  \mathbb{F}\right)  \equiv\mathrm{span}\left\{
1_{\mathbb{D}},i\right\}  \quad\mathbb{C}_{r}\left(  \mathbb{F}\right)
\equiv\mathrm{span}\left\{  1_{\mathbb{D}},j\right\}  \quad\mathbb{C}%
_{g}\left(  \mathbb{F}\right)  \equiv\mathrm{span}\left\{  1_{\mathbb{D}%
},k\right\}  .
\]
We will call these $\ast$-subalgebras the \textbf{blue complex numbers},
\textbf{red complex numbers}, and \textbf{green complex numbers} respectively.
A generic complex number will often be denoted by $z$. The adjugation
operation corresponds to conjugation in the traditional sense over the complex
numbers: if $z=x1_{\mathbb{D}}+yh\in\mathbb{C}_{c}\left(  \mathbb{F}\right)  $
where $h$ is one of $i$, $j$, or $k$, then%
\begin{equation}
z^{\ast}=\left(  x1_{\mathbb{D}}+yh\right)  ^{\ast}=x1_{\mathbb{D}}-yh.
\label{Conjugation generic complex number}%
\end{equation}

Attempts to extend the notion of complex numbers are of course not new: a red
complex number, also called a \textbf{split-complex number} or a
\textbf{hyperbolic number}, was first referenced in the introduction of
\textit{tessarines} by Cockle \cite{Cockle}, and appear in the
split-biquarternions introduced by Clifford \cite{Clifford}. With our
approach, formal symbols are replaced by specific matrices, which makes
computations more explicit and provides a uniform view of all three blue, red,
and green structures.

\subsection{Geometry over complex numbers}

Now we introduce some geometry into the algebra of dihedrons, and then
restrict this to each of the colored complex numbers. This allows us to study a variety of
geometrical situations that are parallel to familiar facts for complex numbers
uniformly under the umbrella of the dihedron geometry.

Define a symmetric bilinear form $\left\langle \cdot,\cdot\right\rangle $ on
$\mathbb{D=D}\left(  \mathbb{F}_{q}\right)  $ by%
\begin{equation}
\left\langle h_{1},h_{2}\right\rangle \equiv\frac{1}{2}\operatorname*{tr}%
\left(  h_{1}h_{2}^{\ast}\right)  \label{Symmetric bilinear form}%
\end{equation}
with associated quadratic form
\begin{equation}
Q\left(  h\right)  \equiv\left\langle h,h\right\rangle =\frac{1}%
{2}\operatorname*{tr}\left(  \det\left(  h\right)  1_{\mathbb{D}}\right)
=\det\left(  h\right).  \nonumber%
\end{equation}
We call the number $Q\left(  h\right)  $ the \textbf{quadrance} of $h$, and
note that $Q\left(  h\right)  =Q\left(  h^{\ast}\right)  $ for any dihedron
$h.$ A key consequence of the connection with the determinant is that for any
dihedrons $h_{1}$ and $h_{2}$ we have
\begin{equation}
Q\left(  h_{1}h_{2}\right)  =Q\left(  h_{1}\right)  Q\left(  h_{2}\right)  .
\label{Quadrance multiplicativity}%
\end{equation}

Specifically the blue, red, and green geometries are obtained by restricting
the symmetric bilinear form $\left\langle \cdot,\cdot\right\rangle $ and
associated quadratic form $Q$ on $\mathbb{D}$ to the respective $\ast
$-subalgebras $\mathbb{C}_{c}\left(  \mathbb{F}\right)  $ for $c=b$, $r$, and
$g.$ If $z_{1}=x_{1}1_{\mathbb{D}}+y_{1}h$ and $z_{2}=x_{2}1_{\mathbb{D}%
}+y_{2}h$, with $h$ either $i,j$ or $k$ in the three cases, then uniformly we
may write
\[
\left\langle z_{1},z_{2}\right\rangle =\frac{1}{2}\operatorname*{tr}\left(
z_{1}z_{2}^{\ast}\right)  =x_{1}x_{2}+\det\left(  h\right)  y_{1}y_{2}.
\]
So in the blue geometry, $\left\langle z_{1},z_{2}\right\rangle
=x_{1}x_{2}+y_{1}y_{2}$ with quadrance $Q\left(  x1_{\mathbb{D}}+yi\right)
=x^{2}+y^{2},$ while in the red and green geometries, $\left\langle
z_{1},z_{2}\right\rangle =x_{1}x_{2}-y_{1}y_{2}$ with quadrance $Q\left(
x1_{\mathbb{D}}+yj\right)  =Q\left(  x1_{\mathbb{D}}+yk\right)  =x^{2}-y^{2}$.
However in the green geometry we can also write the quadrance more simply as
\[
Q\left(
\begin{array}
[c]{cc}%
t & 0\\
0 & u
\end{array}
\right)  =tu.
\]

These formulas make it clear that the blue geometry can be identified with the
familiar Euclidean geometry, while the red and green geometries correspond to
two Einstein-Minkowski (relativistic) geometries. Physicists often refer to
the green geometry in terms of null coordinates or light-cone coordinates in
special relativity theory, by making a change of coordinate from the
$xy$-coordinate to the $tu$-coordinate where $t=x+y$ and $u=x-y$ (see
\cite{Rindler}, \cite{Schwartz}, or \cite{Schwarz}).

\subsection{Unit circles in chromogeometry}

For $c$ one of $b$, $r$, or $g$, the \textbf{unit circle} $S_{c}\left(
\mathbb{F}_{q}\right)  $ is defined to be the set of complex numbers in
$\mathbb{C}_{c,q}$ of quadrance $1_{q}$, that is,%
\begin{equation}
S_{c}\left(  \mathbb{F}_{q}\right)  \equiv\left\{  z\in\mathbb{C}%
_{c,q}\ \colon\ Q\left(  z\right)  =1_{q}\right\}  .
\nonumber%
\end{equation}

From (\ref{Quadrance multiplicativity}) we see that this set contains the
identity $1_{\mathbb{D}},$ is closed under multiplication, and indeed forms a
commutative group, since inverses are given by adjugates and in each of the
three cases all elements are multiples of the identity and a single other matrix.

We identify the blue and red complex numbers $x1_{\mathbb{D}}+yi$ and
$x1_{\mathbb{D}}+yj$ with the affine point $\left[  x,y\right]  $, and
identify the green complex number $x1_{\mathbb{D}}+yk$ with the affine point
$\left[  x+y,x-y\right]  $. With this identification, the unit circles
$S_{c}\left(  \mathbb{F}_{q}\right)  $ can be realized as subsets of
$\mathbb{A}$:%
\begin{align}
S_{b}\left(  \mathbb{F}_{q}\right)   &  =\left\{  \left[  x,y\right]
\in\mathbb{A\ \colon\ }x^{2}+y^{2}=1\right\} \nonumber\\
S_{r}\left(  \mathbb{F}_{q}\right)   &  =\left\{  \left[  x,y\right]
\in\mathbb{A\ \colon\ }x^{2}-y^{2}=1\right\}\nonumber\\
S_{g}\left(  \mathbb{F}_{q}\right)   &  =\left\{  \left[  x,y\right]
\in\mathbb{A\ \colon\ }xy=1\right\}  .\nonumber
\end{align}
There is also the shorthand notations $S_{b,q}$, $S_{r,q}$, and $S_{g,q}$ if
the underlying field is clear from context, and $S_{c,q}$ to indicate the unit
circle of any color. With this point of view, the blue unit circle is the unit
circle in what is the finite field version of Euclidean geometry. The red and
green unit circles can be regarded as two unit circles in finite field
versions of the Minkowski-Einstein (relativistic) geometry.

While the green unit circle $S_{g,q}$ is easily parametrized algebraically,
the traditional parametrizations for the blue and red unit circles over the
real numbers are often presented via trigonometric and hyperbolic
trigonometric functions. However, since our discussion is purely algebraic, we
prefer to utilize the respective rational parametrizations%
\[
\left[  \frac{1-t^{2}}{1+t^{2}},\frac{2t}{1+t^{2}}\right]  \quad
\text{and}\quad\left[  \frac{1+t^{2}}{1-t^{2}},\frac{2t}{1-t^{2}}\right]
\]
for the blue and red unit circles.

These parametrizations conveniently extend also to the finite field case. This
is a classical argument, going back to at least Fermat in the Euclidean case,
and arguably to Euclid from his classification of Pythagorean triples. Recall
that for a prime power $q=p^{r}$, the \textbf{Jacobi symbol} may be defined by%
\begin{equation}
\left(  \frac{-1}{q}\right)  =\left(  \frac{-1}{p}\right)  ^{r}=\left\{
\begin{array}
[c]{cl}%
1 & \text{if }-1\text{ is a square in }\mathbb{F}_{q}\\
-1 & \text{if }-1\text{ is not a square in }\mathbb{F}_{q}\text{.}%
\end{array}
\right.  \nonumber%
\end{equation}

\begin{lemma}
\label{Lem parametrisation of unit circles}Every point on the green unit
circle $S_{g,q}$ can be written as $\left[  x,y\right]  =\left[
t,t^{-1}\right]  $ for a unique element $t\in\mathbb{F}_{q}^{\times}$. Every
point on the red unit circle $S_{r,q}$ except for $\left[  -1,0\right]  ,$ can
be written as $\left[  x,y\right]  =\left[  \frac{1+t^{2}}{1-t^{2}},\frac
{2t}{1-t^{2}}\right]  $ for a unique $t\neq\pm1.$ Every point on the blue unit
circle $S_{b,q},$ except for $\left[  -1,0\right]  ,$ can be written as
$\left[  x,y\right]  =\left[  \frac{1-t^{2}}{1+t^{2}},\frac{2t}{1+t^{2}%
}\right]  $ for a unique $t$ satisfying $t^{2}\neq-1.$ In particular
\[
\left\vert S_{g,q}\right\vert =q-1\quad\left\vert S_{r,q}\right\vert
=q-1\quad\left\vert S_{b,q}\right\vert =q-\left(  \frac{-1}{q}\right)  .
\]

\end{lemma}

\begin{proof}
The parametrization of $S_{g,q}$ is evident with $\left\vert S_{g,q}%
\right\vert =$ $\left\vert \mathbb{F}_{q}^{\times}\right\vert =q-1$. On
$S_{r,q}$, the identity%
\[
\left(  \frac{1+t^{2}}{1-t^{2}}\right)  ^{2}-\left(  \frac{2t}{1-t^{2}%
}\right)  ^{2}=1
\]
holds for all $t\in\mathbb{F}_{q}$ for which $1-t^{2}\neq0$. The line through
the points $\left[  -1,0\right]  $ and $\left[  0,t\right]  $ given by
$y=tx+t$ intersects the red unit circle $x^{2}-y^{2}=1$ in exactly two points,
one of which is $\left[  -1,0\right]  $, and the other, by solving
$x^{2}-\left(  tx+t\right)  ^{2}=1$ for $x$ is $\left[  \frac{1+t^{2}}%
{1-t^{2}},\frac{2t}{1-t^{2}}\right]  $. Hence every point on $S_{r,q}$ except
$\left[  -1,0\right]  $ corresponds to exactly one $t\in\mathbb{F}%
_{q}-\left\{  \pm1\right\}  $, and so $\left\vert S_{r,q}\right\vert
=1+\left(  q-2\right)  =q-1$.

The parametrization of the blue unit circle is done similarly, with a factor
$1+t^{2}$ in the denominator, so we have to exclude those $t\in\mathbb{F}_{q}$
whose square is $-1$, if any. If $-1$ is not a square then the number of such
points is $q+1$. Otherwise, there are exactly two solutions to the equation
$t^{2}=-1$, so the number of points on $S_{b,q}$ is $1+\left(  q-2\right)
=q-1$.
\end{proof}

There is also additional algebraic structure to $\mathbb{C}_{b,q}$ that holds
only when $-1$ is not a square. In this case $\mathbb{C}_{b,q}$ is a field of
size $q^{2}$, since an element $z=x1_{\mathbb{D}}+yi \in \mathbb{C}_{b,q}$ is non-invertible precisely when $Q\left(  z\right)  =x^{2}+y^{2}=\det\left(
z\right)  =0$, and the only solution is $x=y=0$, since if either one of them is
non-zero, say $y$, then $y^{2}\neq0$ and therefore $\left(  \frac{x}%
{y}\right)  ^{2}+1=0$ which is impossible.

\begin{proposition}
\label{Prop cyclic group structure of the unit circles}If $-1$ is a square in
$\mathbb{F}_{q}$ then $S_{g,q}$, $S_{r,q}$, and $S_{b,q}$ are each cyclic
groups of order $q-1$, hence isomorphic. When $-1$ is not a square in
$\mathbb{F}_{q}$, $S_{b,q}$ is a cyclic group of order $q+1$.
\end{proposition}

\begin{proof}
By Lemma \ref{Lem parametrisation of unit circles}, any point on $S_{g,q}$ can
be written as $\left[  t,t^{-1}\right]  $ for $t\in$ $\mathbb{F}_{q}^{\times}%
$. The map $t\mapsto\left[  t,t^{-1}\right]  $ from $\mathbb{F}_{q}^{\times}$
to $S_{g,q}$ is an isomorphism and since $\mathbb{F}_{q}^{\times}$ is cyclic
of order $q-1$, so is $S_{g,q}$. Now $S_{g,q}$ and $S_{r,q}$ are isomorphic
through the map
\[
\theta\ \colon\ S_{r}\left(  \mathbb{F}_{q}\right)  \longrightarrow
S_{g}\left(  \mathbb{F}_{q}\right)  \qquad\mathrm{where}\qquad x1_{\mathbb{D}%
}+yj\mapsto x1_{\mathbb{D}}+yk.
\]
When $-1=\imath^{2}$ for some $\imath\in\mathbb{F}_{q}$, the isomorphism
between $S_{r,q}$ and $S_{b,q}$ is established through the map%
\[
\varrho\ \colon\ S_{b}\left(  \mathbb{F}_{q}\right)  \longrightarrow
S_{r}\left(  \mathbb{F}_{q}\right)  \qquad\mathrm{where}\qquad x1_{\mathbb{D}%
}+yi\mapsto x1_{\mathbb{D}}+y\imath j.
\]
Therefore the three unit circles are cyclic groups of order $q-1$.

Now if $-1$ is not a square in $\mathbb{F}_{q}$, we may extend $\mathbb{F}%
_{q}$ to the quadratic field extension $\mathbb{C}_{b,q}$ and regard $S_{b,q}$
as a subgroup of $\mathbb{C}_{b,q}^{\times}$, hence also cyclic. The order of
$S_{b,q}$ is established by Lemma \ref{Lem parametrisation of unit circles}.
\end{proof}

This cyclic group structure of all of the unit circles has a particularly
pleasant property that will become important later on. We shall state the
result more generally.

\begin{proposition}
\label{Thm cyclic group geom series}Let $\mathbb{F}_{q}$ be a finite field of
size $q$ and $\mathbb{F}_{q}^{\times}$ be its multiplicative group. Suppose
$C$ is a cyclic subgroup of $\mathbb{F}_{q}^{\times}$ of order $n$. Then%
\[
\sum_{t\in C}t^{k}=\left\{
\begin{array}
[c]{cl}%
n1_{q} & \text{if }n~|~k\\
0 & \text{otherwise.}%
\end{array}
\right.
\]

\end{proposition}

\begin{proof}
Let $C$ be generated by $\zeta\neq1_{q}$, so $C=\left\{  1_{q},\zeta,\zeta
^{2},\ldots,\zeta^{n-1}\right\}  $ with $\zeta^{n}=1_{q}$. If $k=mn$ for some
integer $m$, then%
\begin{align*}
\sum_{t\in C}t^{k}  &  =1_{q}^{mn}+\zeta^{mn}+\zeta^{2mn}+\cdots
+\zeta^{\left(  n-1\right)  mn}\\
&  =\underset{n\text{ times}}{\underbrace{1_{q}+1_{q}+\cdots+1_{q}}}=n1_{q}.
\end{align*}

If $n$ does not divide $k$, then $\zeta^{k}\neq1_{q}$ and so $\displaystyle\sum_{t\in C}t^{k}=\sum_{l=0}^{n-1}\zeta^{kl}=\frac{\zeta^{kn}-1_{q}}%
{\zeta^{k}-1_{q}}=0$ since $\zeta^{kn}=\left(  \zeta^{n}\right)  ^{k}=1_{q}$ and $\zeta^{k}%
\neq1_{q}$.
\end{proof}

\subsection{Special orthogonal groups in chromogeometry}

The three different geometries on $\mathbb{C}_{c,q}$ induced by the symmetric
bilinear form (\ref{Symmetric bilinear form}) correspond to the green, red,
and blue geometries on $\mathbb{A}$ by the identification stipulated
previously: $x1_{\mathbb{D}}+yi$ and $x1_{\mathbb{D}}+yj$ as $\left[
x,y\right]  $, and $x1_{\mathbb{D}}+yk$ as $\left[  x+y,x-y\right]  $.

More specifically, over the green geometry, for $z_{1}=x_{1}1_{\mathbb{D}%
}+y_{1}k$ and $z_{2}=x_{2}1_{\mathbb{D}}+y_{2}k\in\mathbb{C}_{g,q}$ we have
\[
\left\langle z_{1},z_{2}\right\rangle =\frac{1}{2}\operatorname*{tr}\left(
z_{1}z_{2}^{\ast}\right)  =\frac{1}{2}\left(  \left(  x_{1}+y_{1}\right)
\left(  x_{2}-y_{2}\right)  +\left(  x_{1}-y_{1}\right)  \left(  x_{2}%
+y_{2}\right)  \right)
\]
so the corresponding symmetric bilinear form on $\mathbb{A}$ is%
\[
\left[  x_{1},y_{1}\right]  \cdot_{g}\left[  x_{2},y_{2}\right]  \equiv
\frac{1}{2}\left(  x_{1}y_{2}+x_{2}y_{1}\right)  .
\]
We call the geometry induced by $\cdot_{g}$ the green geometry on $\mathbb{A}$.

Similarly, we have the red and blue geometry on $\mathbb{A}$ defined by the
symmetric bilinear forms%
\begin{align*}
\left[  x_{1},y_{1}\right]  \cdot_{r}\left[  x_{2},y_{2}\right]   &  \equiv
x_{1}x_{2}-y_{1}y_{2}\\
\left[  x_{1},y_{1}\right]  \cdot_{b}\left[  x_{2},y_{2}\right]   &  \equiv
x_{1}x_{2}+y_{1}y_{2}.
\end{align*}
A simple algebraic manipulation shows that $\left[  x_{1},y_{1}\right]  \cdot_{c}\left[  x_{2},y_{2}\right]  =\left(
x_{1},y_{1}\right)  M_{c}\left(  x_{2},y_{2}\right)^{T}$ where $M_{c}$ are the invertible matrices%
\begin{equation}
M_{g}=\left(
\begin{array}
[c]{cc}%
0 & \frac{1}{2}\\
\frac{1}{2} & 0
\end{array}
\right)  \quad M_{r}=\left(
\begin{array}
[c]{cc}%
1 & 0\\
0 & -1
\end{array}
\right)  \quad M_{b}=\left(
\begin{array}
[c]{cc}%
1 & 0\\
0 & 1
\end{array}
\right)  . \label{List of Mc}%
\end{equation}

If $\mathrm{GL}(2,\mathbb{F}_{q})$ denotes the group of $2\times2$ invertible
matrices over $\mathbb{F}_{q}$, there is a natural right action of
$\mathrm{GL}(2,\mathbb{F}_{q})$ on $\mathbb{A}$ and left action of
$\mathrm{GL}(2,\mathbb{F}_{q})$ on $\mathrm{Pol}_{2}\left(  \mathbb{F}%
_{q}\right)  $ as%
\begin{align}
\left[  x,y\right]  \cdot\left(
\begin{array}
[c]{cc}%
h_{11} & h_{12}\\
h_{21} & h_{22}%
\end{array}
\right)   &  \equiv\left[  h_{11}x+h_{21}y,h_{12}x+h_{22}y\right] \nonumber\\
\left(
\begin{array}
[c]{cc}%
h_{11} & h_{12}\\
h_{21} & h_{22}%
\end{array}
\right)  \cdot\pi &  \equiv\pi\left(  h_{11}\alpha+h_{21}\beta,h_{12}%
\alpha+h_{22}\beta\right)\nonumber
\end{align}
where we regard $\pi$ as a polynumber in $\alpha$ and $\beta$. The subgroup of
$\mathrm{GL}(2,\mathbb{F}_{q})$ consisting of matrices $h$ that preserve that
bilinear form $\cdot_{c}$ is called the \textbf{orthogonal group}
$O_{c}\left(  \mathbb{F}_{q}\right)  =O_{c,q}$ with respect to geometry $c$.
Any matrix $h\in O_{c,q}$ necessarily satisfies $hM_{c}h^{T}=M_{c}$ where
$M_{c}$ is given in (\ref{List of Mc}) and therefore $\det\left(  h\right)
=\pm1$.

The collection of matrices with determinant $1$ in $O_{c,q}$ is itself a group
which we call the \textbf{rotation group} $SO_{c}\left(  \mathbb{F}%
_{q}\right)  =SO_{c,q}$. We are particularly interested in the analysis
pertaining to the rotation group $SO_{c,q}$ by considering how $SO_{c,q}$ acts
on various spaces. $SO_{c,q}$ naturally right acts on $\mathbb{A}$ via the
restriction of the right action of $\mathrm{GL}(2,\mathbb{F}_{q})$ on
$\mathbb{A}$. There is then an associated left action of $SO_{c,q}$ on
$\mathbb{F}_{q}^{\mathbb{A}}$ and an associated right action of $SO_{c,q}$ on
$\left(  \mathbb{F}_{q}^{\mathbb{A}}\right)  ^{\ast}$. Similarly, the left
action of $SO_{c,q}$ on $\mathrm{Pol}_{2}\left(  \mathbb{F}_{q}\right)  $ is
induced as the restriction of the action of $\mathrm{GL}(2,\mathbb{F}_{q})$ on
$\mathrm{Pol}_{2}\left(  \mathbb{F}_{q}\right)  .$ This action respects
evaluation in a sense that%
\begin{equation}
\varepsilon\left(  h\cdot\pi\right)  =h\cdot\varepsilon\left(  \pi\right)  .
\label{SO action respects evaluation}%
\end{equation}
This setup now makes it easy and pleasant to explicitly describe the rotation
group $SO_{c,q}$.

\begin{theorem}
\label{Thm special orthogonal groups in three geometries}The rotation group
$SO_{c,q}$ is precisely the unit circle group $S_{c,q}$.
\end{theorem}

\begin{proof}
Any $h\in SO_{c,q}$ satisfies $\det\left(  h\right)  =1$ by definition. Since
over any geometry $c$, $h$ must be in $\mathbb{C}_{c,q}$, the condition
$\det\left(  h\right)  =1$ is equivalent to $h\in S_{c,q}$. This implies that
$SO_{c,q}\subseteq S_{c,q}$.

Similarly, if $z\in S_{c,q}$ then $Q\left(  z\right)  =\det\left(  z\right)
=1$ by definition. The condition that $z$ preserves the bilinear form
$\cdot_{c}$ on $\mathbb{A}$ is equivalent to $zM_{c}z^{T}=M_{c}$. Over the
green geometry, $z$ has the form $x1_{\mathbb{D}}+yk$ where $x^{2}-y^{2}=1$
and so by inspection,%
\[
zM_{g}z^{T}=\left(
\begin{array}
[c]{cc}%
0 & \frac{1}{2}\left(  x^{2}-y^{2}\right) \\
\frac{1}{2}\left(  x^{2}-y^{2}\right)  & 0
\end{array}
\right)  =\left(
\begin{array}
[c]{cc}%
0 & \frac{1}{2}\\
\frac{1}{2} & 0
\end{array}
\right)  =M_{g}.
\]
The cases $z\in S_{r,q}$ and $z\in S_{b,q}$ are treated similarly. In all
cases, $zM_{c}z^{T}=M_{c}$ and thus $z$ preserves the bilinear form $\cdot
_{c}$ on $\mathbb{A}$. This shows that $S_{c,q}\subseteq SO_{c,q}$.
\end{proof}

The result above allows us to present the following observation about the
action of $SO_{c,q}$ on $\mathbb{A}$.

\begin{theorem}
\label{S is an SO orbit}The unit circle $S_{c,q}$ is an orbit of the group
$SO_{c,q}$ in $\mathbb{A}$. Moreover, the action of $SO_{c,q}$ in $S_{c,q}$ is faithful.
\end{theorem}

\begin{proof}
Since $S_{c,q}=SO_{c,q}$ by Theorem
\ref{Thm special orthogonal groups in three geometries}, we may identify
$\left[  x,y\right]  \cdot h$ as $zh$ where $z\in\mathbb{C}_{c,q}$ is the
complex number associated to $\left[  x,y\right]  $. This allows us to observe
that $\left[  x,y\right]  \cdot h\in S_{c,q}$ for any $h\in SO_{c,q}$. To see
that the action of $SO_{c,q}$ in $S_{c,q}$ is faithful, the equation $\left[
x,y\right]  \cdot h=\left[  x,y\right]  $ for any $\left[  x,y\right]  \in
S_{c,q}$ may be identified as $zh=z$ for all $z\in S_{c,q}$. This yields
$h=1_{\mathbb{D}}$, finishing our claim.
\end{proof}

As a final remark, the theory of orthogonal matrices over finite fields over a
general geometry has been intensively studied (see \cite{Carlitz},
\cite{MacWilliams2}, \cite{Taylor}). The size of the group $SO_{b}\left(
n,\mathbb{F}_{q}\right)  $ consisting of all $n\times n$ orthogonal matrices
over $\mathbb{F}_{q}$ with determinant $1$ is given by
\[
\left\vert SO_{b}\left(  2m+1,\mathbb{F}_{q}\right)  \right\vert =q^{m^{2}}%
{\displaystyle\prod\limits_{k=1}^{m}}
\left(  q^{2k}-1\right)
\]
and%
\[
\left\vert SO_{b}\left(  2m,\mathbb{F}_{q}\right)  \right\vert =\left\{
\begin{array}
[c]{cl}%
\dfrac{1}{q^{m}+1}%
{\displaystyle\prod\limits_{k=0}^{m-1}}
\left(  q^{2m}-q^{2k}\right)  & \text{if }\left(  \frac{-1}{q}\right)  =1\\
\dfrac{1}{q^{m}+\left(  -1\right)  ^{m}}%
{\displaystyle\prod\limits_{k=0}^{m-1}}
\left(  q^{2m}-q^{2k}\right)  & \text{if }\left(  \frac{-1}{q}\right)
=-1\text{.}%
\end{array}
\right.
\]
It is clear that when $n=2$, the formula above matches with the result
established in Theorem \ref{Thm special orthogonal groups in three geometries}.

\section{Circular Integral Functionals}

We now introduce an algebraic integral on the unit circles in each of the
three geometries, valid over a general finite field $\mathbb{F}_{q}$ of odd
characteristic. These will be functionals on polynumbers which satisfies
conditions that we call \textit{normalization}, \textit{locality }and
\textit{invariance}. We begin with the simpler case of a corresponding
functional on the space $\mathbb{F}_{q}^{\mathbb{A}}$ of $\mathbb{F}_q$-valued
functions on $\mathbb{A}$.

\subsection{Circular integral functionals on $\mathbb{F}_{q}^{\mathbb{A}}$}

For one of the three geometries $c=b$, $r$, or $g$ on the affine space
$\mathbb{A}=\mathbb{A}\left(  \mathbb{F}_{q}\right)  $, the unit circle is
$S_{c,q}$ and the special orthogonal group is $SO_{c,q}$. A linear functional
$\phi$ on $\mathbb{F}_{q}^{\mathbb{A}}$ is a \textbf{circular integral
functional }precisely when the following three conditions are satisfied:

\begin{description}
\item[(Normalization)] The functional $\phi$ maps the constant function
$1_{\mathbb{F}^{\mathbb{A}}}$ to $1_{q}$

\item[(Locality)] If $f\in\mathbb{F}_{q}^{\mathbb{A}}$ restricts to the zero
function on $S_{c,q}$, then $\phi\left(  f\right)  =0$

\item[(Invariance)] The functional $\phi$ is $SO_{c,q}$-invariant, that is
$\phi\cdot h=\phi$ for all $h\in SO_{c,q}$.
\end{description}

\begin{theorem}
\label{Existence and Uniqueness of CIF on F^A}There is one and only one
circular integral functional on $\mathbb{F}_{q}^{\mathbb{A}}$, namely%
\[
\phi\left(  f\right)  =\frac{1}{\left\vert S_{c,q}\right\vert }\sum_{\left[
x,y\right]  \in S_{c,q}}f\left(  x,y\right)  .
\]

\end{theorem}

\begin{proof}
The definition of $\phi$ above is well-defined since $\left\vert
S_{c,q}\right\vert \not \equiv 0~\operatorname{mod}p$ by Proposition \ref{Prop cyclic group structure of the unit circles}. This functional is
normalized, since
\[
\phi\left(  1_{\mathbb{F}^{\mathbb{A}}}\right)  =\frac{1}{\left\vert
S_{c,q}\right\vert }\sum_{\left[  x,y\right]  \in S_{c,q}}1_{q}=1_{q}%
\]
and from the definition it is clearly localized. For the Invariance property,
for $h\in SO_{c,q}$ and $f\in\mathbb{F}_{q}^{\mathbb{A}}$, we have $\left(
\phi\cdot h\right)  \left(  f\right)  =\phi\left(  h\cdot f\right)  $ and
therefore
\[
\left(  \phi\cdot h\right)  \left(  f\right)  =\frac{1}{\left\vert
S_{c,q}\right\vert }\sum_{\left[  x,y\right]  \in S_{c,q}}\left(  h\cdot
f\right)  \left(  x,y\right)  =\frac{1}{\left\vert S_{c,q}\right\vert }%
\sum_{\left[  x,y\right]  \in S_{c,q}}f\left(  \left[  x,y\right]  \cdot
h\right)  =\frac{1}{\left\vert S_{c,q}\right\vert }\sum_{\left[  x,y\right]
\in S_{c,q}}f\left(  x,y\right)  =\phi\left(  f\right)
\]
since $SO_{c,q}$ acts as a group of bijections on $S_{c,q}$. Consequently
$\phi$ is a circular integral functional.

Now if $\eta$ is any circular integral functional on $\mathbb{F}%
_{q}^{\mathbb{A}},$ since $\Delta$ is a basis for $\mathbb{F}_{q}^{\mathbb{A}%
}$, it suffices to show that $\eta\left(  \delta_{\left[  a,b\right]
}\right)  =\phi\left(  \delta_{\left[  a,b\right]  }\right)  $ for all
$\left[  a,b\right]  \in\mathbb{A}$. We observe that
\begin{equation}
\phi\left(  \delta_{\left[  a,b\right]  }\right)  =\frac{1}{\left\vert
S_{c,q}\right\vert }\sum_{\left[  x,y\right]  \in S_{c,q}}\delta_{\left[
a,b\right]  }\left(  x,y\right)  =\left\{
\begin{array}
[c]{cl}%
\frac{1}{\left\vert S_{c,q}\right\vert } & \text{if }\left[  a,b\right]  \in
S_{c,q}\\
0 & \text{if }\left[  a,b\right]  \notin S_{c,q}.
\end{array}
\right.  \label{A formula for phi at delta function}%
\end{equation}

If $\left[  a,b\right]  \notin S_{c,q}$, then $\delta_{\left[  a,b\right]
}=0$ on $S_{c,q}$ and therefore by Locality, $\eta\left(  \delta_{\left[
a,b\right]  }\right)  =0$. On the other hand if $\left[  a,b\right]  \in
S_{c,q}$, then by Invariance, for any $h\in SO_{c,q}$ it follows that%
\begin{equation}
\eta\left(  \delta_{\left[  a,b\right]  }\right)  =\eta\left(  h\cdot
\delta_{\left[  a,b\right]  }\right)  =\eta\left(  \delta_{\left[  a,b\right]
\cdot h}\right)  . \label{G-invariance on eta}%
\end{equation}
Since $SO_{c,q}$ acts transitively on $S_{c,q}$, (\ref{G-invariance on eta})
implies that $\eta$ is constant on the set $\left.  \Delta\right\vert
_{S_{c,q}}$. Consequently,
\[
1_{q}=\eta\left(  1_{\mathbb{F}^{\mathbb{A}}}\right)  =\eta\left(
\sum_{\left[  x,y\right]  \in\mathbb{A}}\delta_{\left[  x,y\right]  }\right)
=\eta\left(  \sum_{\left[  x,y\right]  \in S_{c,q}}\delta_{\left[  x,y\right]
}\right)  =\sum_{\left[  x,y\right]  \in S_{c,q}}\eta\left(  \delta_{\left[
x,y\right]  }\right)  =\left\vert S_{c,q}\right\vert \eta\left(
\delta_{\left[  a,b\right]  }\right)
\]
where the second equality comes from (\ref{Property of constant 1 function}).
This concludes the proof.
\end{proof}

\subsection{Circular integral functionals on $\mathrm{Pol}_{2}\left(
\mathbb{F}_{q}\right)  $}

A linear functional $\psi$ on $\mathrm{Pol}_{2}\left(  \mathbb{F}_{q}\right)
$ is called a \textbf{circular integral functional }precisely when the
following three conditions are satisfied:

\begin{description}
\item[(Normalization)] The functional $\psi$ maps the polynumber $\mathbf{1}$
to $1_{q}$.

\item[(Locality)] If $\pi\in\mathrm{Pol}_{2}\left(  \mathbb{F}_{q}\right)  $
evaluates to the zero function on $S_{c,q}$ then $\psi\left(  \pi\right)  =0$.

\item[(Invariance)] The functional $\psi$ is $SO_{c,q}$-invariant, that is
$\psi\left(  h\cdot\pi\right)  =\psi\left(  \pi\right)  $ for any $\pi
\in\mathrm{Pol}_{2}\left(  \mathbb{F}_{q}\right)  $ and $h\in SO_{c,q}$.
\end{description}

These three conditions are virtually identical to the corresponding conditions
for a circular integral functional on $\mathbb{F}_{q}^{\mathbb{A}}$, however
there is quite a subtle difference in the analysis. Recall that $\varepsilon\left(  \pi\right)  $ refers to the evaluation of a polynumber
$\pi.$

\begin{theorem}
\label{Existence and Uniqueness of CIF on Pol space}There is one and only one
circular integral functional on $\mathrm{Pol}_{2}\left(  \mathbb{F}%
_{q}\right)  $ namely%
\begin{equation}
\psi_{c,q}\left(  \pi\right)  =\frac{1}{\left\vert S_{c,q}\right\vert }%
\sum_{\left[  x,y\right]  \in S_{c,q}}\varepsilon\left(  \pi\right)  \left(
x,y\right)  . \nonumber%
\end{equation}

\end{theorem}

\begin{proof}
We first note that $\psi_{c,q}=\phi\circ\varepsilon$ where $\phi$ is the
unique circular integral functional on $\mathbb{F}_{q}^{\mathbb{A}}$
established in Theorem \ref{Existence and Uniqueness of CIF on F^A}. Since
$\varepsilon\left(  \mathbf{1}\right)  =1_{\mathbb{F}^{\mathbb{A}}}$, the
Normalization condition follows from that of $\phi$. The Locality property is
also obvious, for if $\pi\in\mathrm{Pol}_{2}\left(  \mathbb{F}_{q}\right)  $
evaluates to the zero function on $S_{c,q}$, then $\varepsilon\left(
\pi\right)  \left(  x,y\right)  =0$ for all $\left[  x,y\right]  \in S_{c,q}$.

To prove the Invariance property, we note from (\ref{SO action respects evaluation}) that the $SO_{c,q}$-left-action on
$\mathrm{Pol}_{2}\left(  \mathbb{F}_{q}\right)  $ respects evaluation, so%
\begin{align*}
\psi_{c,q}\left(  h\cdot\pi\right)   &  =\frac{1}{\left\vert S_{c,q}%
\right\vert }\sum_{\left[  x,y\right]  \in S_{c,q}}\varepsilon\left(
h\cdot\pi\right)  \left(  x,y\right)  =\frac{1}{\left\vert S_{c,q}\right\vert
}\sum_{\left[  x,y\right]  \in S_{c,q}}\left(  h\cdot\varepsilon\left(
\pi\right)  \right)  \left(  x,y\right) \\
&  =\phi\left(  h\cdot\varepsilon\left(  \pi\right)  \right)  =\phi\left(
\varepsilon\left(  \pi\right)  \right)  =\psi_{c,q}\left(  \pi\right)
\end{align*}
where we used the Invariance property of $\phi$. We deduce that $\psi_{c,q}$
is a circular integral functional on $\mathrm{Pol}_{2}\left(  \mathbb{F}%
_{q}\right)  $ with respect to $S_{c,q}$.

Now if $\chi$ is any circular integral functional on $\mathrm{Pol}_{2}\left(
\mathbb{F}_{q}\right)  $ with respect to $S_{c,q}$, we need to show that
$\chi\left(  \pi\right)  =\psi_{c,q}\left(  \pi\right)  $ for any $\pi\in$
$\mathrm{Pol}_{2}\left(  \mathbb{F}_{q}\right)  $. By Theorem
\ref{Decomposition of Pol space}, the decomposition $\pi=\pi_{1}+\pi_{2}$
where $\pi_{1}\in\mathcal{P}_{q}$ and $\pi_{2}\in\ker\left(  \varepsilon
\right)  $ is unique. Since $\varepsilon\left(  \pi_{2}\right)  =0$, Locality
asserts that $\chi\left(  \pi_{2}\right)  =0$. This means that we only need to
consider the possible values of $\chi$ on $\mathcal{P}_{q}$.

Since $\Pi$ is a basis for $\mathcal{P}_{q}$, it suffices to show that%
\[
\chi\left(  \pi_{\left[  a,b\right]  }\right)  =\psi_{c,q}\left(  \pi_{\left[
a,b\right]  }\right)  =\left\{
\begin{array}
[c]{cl}%
\frac{1}{\left\vert S_{c,q}\right\vert } & \text{if }\left[  a,b\right]  \in
S_{c,q}\\
0 & \text{if }\left[  a,b\right]  \notin S_{c,q}%
\end{array}
\right.
\]
where the second equality follows from (\ref{A formula for phi at delta function}).

If $\left[  a,b\right]  \notin S_{c,q}$, then $\varepsilon\left(  \pi_{\left[
a,b\right]  }\right)  =\delta_{\left[  a,b\right]  }=0$ on $S_{c,q}$ which
implies that $\chi\left(  \pi_{\left[  a,b\right]  }\right)  =0$ by Locality.
Now suppose $\left[  a,b\right]  \in S_{c,q}$. In this case,%
\[
\varepsilon\left(  h\cdot\pi_{\left[  a,b\right]  }\right)  =h\cdot
\varepsilon\left(  \pi_{\left[  a,b\right]  }\right)  =h\cdot\delta_{\left[
a,b\right]  }=\delta_{\left[  a,b\right]  \cdot h^{-1}}=\varepsilon\left(
\pi_{\left[  a,b\right]  \cdot h^{-1}}\right)
\]
where the third equality above follows from equation
(\ref{Action of G to delta function}). It follows that $h\cdot\pi_{\left[
a,b\right]  }-\pi_{\left[  a,b\right]  \cdot h^{-1}}\in\ker\left(
\varepsilon\right)  $. By Invariance and Locality,%
\[
\chi\left(  \pi_{\left[  a,b\right]  \cdot h^{-1}}\right)  =\chi\left(
h\cdot\pi_{\left[  a,b\right]  }\right)  =\chi\left(  \pi_{\left[  a,b\right]
}\right)  .
\]
This implies that $\chi$ is constant on $\left.  \Pi\right\vert _{S_{c,q}}$.

Finally,
\[
1_{q}=\chi\left(  \mathbf{1}\right)  =\chi\left(  \sum_{\left[  a,b\right]
\in\mathbb{A}}\pi_{\left[  a,b\right]  }\right)  =\chi\left(  \sum_{\left[
a,b\right]  \in S_{c,q}}\pi_{\left[  a,b\right]  }\right)  =\sum_{\left[
a,b\right]  \in S_{c,q}}\chi\left(  \pi_{\left[  a,b\right]  }\right)
=\left\vert S_{c,q}\right\vert \chi\left(  \pi_{\left[  a,b\right]  }\right)
\]
which concludes the proof.
\end{proof}

We call the circular integral functionals $\psi_{c,q}$ the \textbf{green}, \textbf{red}, or \textbf{blue Fourier summation functionals} respectively.
Since $\left\vert S_{c,q}\right\vert $ is known, we obtain the following three
formulas:%

\begin{align*}
\psi_{g,q}\left(\alpha^{k}\beta^{l}\right) \equiv-\sum\limits_{\left[  x,y\right]  \in S_{g,q}} x^{k}y^{l} \qquad \psi_{r,q}\left(\alpha^{k}\beta^{l}\right) \equiv-\sum_{\left[x,y\right] \in S_{r,q}}x^{k}y^{l} \qquad \psi_{b,q}\left(\alpha^{k}\beta^{l}\right) \equiv -\left(\frac{-1}{q}\right) \sum_{\left[x,y\right] \in S_{b,q}}x^{k}y^{l}.
\end{align*}%

We give two examples of the computation of the Fourier summation functionals
$\psi_{c,q}$.

\begin{example}
\label{Example q = 17}Over $\mathbb{F}_{17}=\left\{  0,1,2,\ldots
,15,16\right\}  $, the red unit circle is%
\[
S_{r,17}=\left\{
\begin{array}
[c]{c}%
\left[  1,0\right]  \ \left[  16,0\right]  \ \left[  6,1\right]  \ \left[
6,16\right]  \ \left[  11,1\right]  \ \left[  11,16\right]  \ \left[
3,5\right]  \ \left[  3,12\right] \\
\left[  14,5\right]  \ \left[  14,12\right]  \ \left[  4,7\right]  \ \left[
4,10\right]  \ \left[  13,7\right]  \ \left[  13,10\right]  \ \left[
0,4\right]  \ \left[  0,13\right]
\end{array}
\right\}  .
\]

It can be verified that, for example,%
\begin{align*}
\psi_{r,17}\left(  \alpha^{6}\beta^{4}\right)   &  =-\sum_{\left[  x,y\right]
\in S_{r,17}}x^{6}y^{4}\\
&  =-\left(
\begin{array}
[c]{c}%
1^{6}0^{4}+16^{6}0^{4}+0^{6}1^{4}+0^{6}16^{4}+4^{6}6^{4}+4^{6}11^{4}%
+13^{6}6^{4}+13^{6}11^{4}+\\
6^{6}4^{4}+6^{6}13^{4}+11^{6}4^{4}+11^{6}13^{4}+3^{6}3^{4}+3^{6}14^{4}%
+14^{6}3^{4}+14^{6}14^{4}%
\end{array}
\right)  =3.
\end{align*}

\end{example}

\begin{example}
\label{Example q = 27}Since $1+\alpha-\alpha^{3}$ is irreducible in
$\mathrm{Pol}\left(  \mathbb{F}_{3}\right)  $, we may regard $\mathbb{F}_{27}$
as $\mathrm{Pol}\left(  \mathbb{F}_{3}\right)  /\left(  1+\alpha-\alpha
^{3}\right)  $. By identifying $\pi+\left(  1+\alpha-\alpha^{3}\right)  $
simply as $\pi$, the blue unit circle is
\begingroup
\allowdisplaybreaks
\[
S_{b,27}=\left\{
\begin{array}
[c]{c}%
\left[  1,0\right]  \ \left[  2,0\right]  \ \left[  0,1\right]  \ \left[
0,2\right]  \ \left[  \alpha^{2},2\alpha+\alpha^{2}\right]  \ \left[
\alpha^{2},\alpha+2\alpha^{2}\right]  \ \left[  2\alpha^{2},2\alpha+\alpha
^{2}\right] \\
\ \left[  2\alpha^{2},\alpha+2\alpha^{2}\right]  \ \left[  2\alpha+\alpha
^{2},\alpha^{2}\right]  \ \left[  \alpha+2\alpha^{2},\alpha^{2}\right]
\ \left[  2\alpha+\alpha^{2},2\alpha^{2}\right] \\
\left[  \alpha+2\alpha^{2},2\alpha^{2}\right]  \ \left[  2+\alpha^{2}%
,1+\alpha+\alpha^{2}\right]  \ \left[  2+\alpha^{2},2+2\alpha+2\alpha
^{2}\right] \\
\ \left[  2\alpha+2\alpha^{2},1+2\alpha+\alpha^{2}\right]  \ \ \left[
1+2\alpha^{2},2+2\alpha+2\alpha^{2}\right] \\
\left[  1+2\alpha^{2},1+\alpha+\alpha^{2}\right]  \ \left[  1+\alpha
+\alpha^{2},2+\alpha^{2}\right]  \ \left[  2+2\alpha+2\alpha^{2},2+\alpha
^{2}\right] \\
\left[  1+\alpha+\alpha^{2},1+2\alpha^{2}\right]  \ \left[  \alpha+\alpha
^{2},1+2\alpha+\alpha^{2}\right]  \ \left[  \alpha+\alpha^{2},2+\alpha
+2\alpha^{2}\right] \\
\ \left[  2\alpha+2\alpha^{2},2+\alpha+2\alpha^{2}\right]  \ \left[
1+2\alpha+\alpha^{2},\alpha+\alpha^{2}\right]  \ \left[  2+\alpha+2\alpha
^{2},\alpha+\alpha^{2}\right] \\
\ \left[  2+2\alpha+2\alpha^{2},1+2\alpha^{2}\right]  \ \left[  1+2\alpha
+\alpha^{2},2\alpha+2\alpha^{2}\right]  \ \left[  2+\alpha+2\alpha^{2}%
,2\alpha+2\alpha^{2}\right]
\end{array}
\right\}  .
\]
\endgroup
Using this we may directly compute that%
\[
\psi_{b,27}\left(  \alpha^{4}\beta^{6}\right)  =-\left(  \frac{-1}{27}\right)
\sum_{\left[  x,y\right]  \in S_{b,27}}x^{4}y^{6}=-\left(  \frac{-1}%
{3}\right)  ^{3}\sum_{\left[  x,y\right]  \in S_{b,27}}x^{4}y^{6}=0.
\]

\end{example}

\section{Explicit Formulas for $\psi_{c,q}$}

The computation of $\psi_{c,q}$ can become imposing, especially when the size
of the field is large. The second half of this paper develops explicit
formulas that turn out to involve the circular super Catalan numbers in the
red and blue geometries. A formula for the green Fourier summation functional
$\psi_{g,q}$ is relatively easy to find, and will be the key to obtaining
formulas for $\psi_{r,q},$ and for $\psi_{b,q}$ in the case when $-1$ is a
square in $\mathbb{F}_{q}$. This chain of connection suggests that
relativistic geometries might be simpler than Euclidean geometry, at least in
the analysis context of integration theories over general fields.

The formula for $\psi_{b,q}$ when $-1$ is not a square in $\mathbb{F}_{q}$
warrants a special treatment involving Fourier analysis style arguments that
are closer to the case for the field of real numbers, in which $-1$ is also
not a square.

\begin{theorem}
\label{Thm CIF green group structure}In the green geometry over $\mathbb{F}%
_{q}$, for any $k,l\in\mathbb{N}$,%
\begin{equation}
\psi_{g,q}\left(  \alpha^{k}\beta^{l}\right)  =\left\{
\begin{array}
[c]{cl}%
1_{q} & \text{if }\left\vert S_{g,q}\right\vert \mid\left(  k-l\right)  ,\\
0 & \text{otherwise.}%
\end{array}
\right.\nonumber
\end{equation}

\end{theorem}

\begin{proof}
From the parametrization of $S_{g,q}$, we observe that
\[
\psi_{g,q}\left(  \alpha^{k}\beta^{l}\right)  =-\sum_{\left[  x,y\right]  \in
S_{g,q}}x^{k}y^{l}=-\sum_{t\in\mathbb{F}_{q}^{\times}}t^{k}\left(
t^{-1}\right)  ^{l}=-\sum_{t\in\mathbb{F}_{q}^{\times}}t^{k-l}.
\]
The conclusion follows by applying Lemma \ref{Thm cyclic group geom series}
with $C=\mathbb{F}_{q}^{\times}$ since $q-1$ equals $\left\vert S_{g,q}%
\right\vert .$
\end{proof}

An explicit formula for $\psi_{r,q}$ can then be found by utilizing the
isomorphism $\theta$ between $S_{r,q}$ and $S_{g,q}$, and involves two
numerical invariants along with the coefficients of the circular polynumbers
$\pi_{k,l}.$

\begin{theorem}
\label{Theorem formula red CIF}\ In the red geometry over $\mathbb{F}_{q}$,
define $w\equiv\frac{1}{2}\left\vert S_{b,q}\right\vert =\frac{1}{2}\left(
q-1\right)  $, and for natural numbers $k$ and $l$, let $R_{k,l}\equiv
\left\lfloor \frac{k+l}{\left\vert S_{b,q}\right\vert }\right\rfloor .$ Then
\begin{equation}
\psi_{r,q}(\alpha^{k}\beta^{l})=\left\{
\begin{array}
[c]{cl}%
\left(
{\displaystyle\sum\limits_{\left\vert d\right\vert \leq R_{k,l}}}
\left[  \alpha^{m+n+dw}\right]  \right)  \pi_{k,l}\operatorname{mod}p &
\text{if }k=2m\text{, }l=2n\\
0 & \text{otherwise.}%
\end{array}
\right.\nonumber
\end{equation}

\end{theorem}

\begin{proof}
For any $\left[  x,y\right]  \in S_{r,q}$, write $u=x+y$ and $v=x-y$ so that
$\left[  u,v\right]  \in S_{g,q}$. Inversely the expressions $x=\frac{u+v}{2}$
and $y=\frac{u-v}{2}$ are both well-defined since $2\neq0$. Using the
definition of $\psi_{r,q}$,
\[
\psi_{r,q}\left(  \alpha^{k}\beta^{l}\right)  =-\sum_{\left[  x,y\right]  \in
S_{r,q}}x^{k}y^{l}=-\frac{1}{2^{k+l}}\sum_{\left[  u,v\right]  \in S_{g,q}%
}\left(  u+v\right)  ^{k}\left(  u-v\right)  ^{l}.
\]
Reverting back to $x$ and $y$, we obtain a formula for $\psi_{r,q}$ in terms
of $\psi_{g,q}$ as follows:
\begin{align}
\psi_{r,q}\left(  \alpha^{k}\beta^{l}\right)   &  =-\frac{1}{2^{k+l}}%
\sum_{\left[  x,y\right]  \in S_{g,q}}\left(  \sum_{s=0}^{k}\binom{k}{s}%
x^{s}y^{k-s}\right)  \left(  \sum_{t=0}^{l}\left(  -1\right)  ^{l-t}\binom
{l}{t}x^{t}y^{l-t}\right) \nonumber\\
&  =-\frac{1}{2^{k+l}}\sum_{s=0}^{k}\sum_{t=0}^{l}\left(  -1\right)
^{l-t}\binom{k}{s}\binom{l}{t}\sum_{\left[  x,y\right]  \in S_{g,q}}%
x^{s+t}y^{k+l-\left(  s+t\right)  }\nonumber\\
&  =\frac{\left(  -1\right)  ^{l}}{2^{k+l}}\sum_{s=0}^{k}\sum_{t=0}^{l}\left(
-1\right)  ^{t}\binom{k}{s}\binom{l}{t}\psi_{g,q}\left(  \alpha^{s+t}%
\beta^{k+l-\left(  s+t\right)  }\right)  \label{A form of psi red part 1}%
\end{align}
where the minus sign originally in front of the expression has been absorbed.

By Theorem \ref{Thm CIF green group structure}, $\psi_{g,q}\left(
\alpha^{s+t}\beta^{k+l-\left(  s+t\right)  }\right)  =1_{q}$ precisely when%
\begin{equation}
2\left(  s+t\right)  =k+l+d\left(  q-1\right)
\label{Divisibility condition psi red}%
\end{equation}
for some $d\in\mathbb{Z}$ and $0$ otherwise. Note that here $s$ and $t$ are
indices in the ranges $0\leq s\leq k$ and $0\leq t\leq l.$

If $k+l$ is odd, then by comparing the parity of both sides, (\ref{Divisibility condition psi red}) has no solution for $s$ and $t$. We
conclude that $\psi_{g,q}\left(  \alpha^{s+t}\beta^{k+l-\left(  s+t\right)
}\right)  =0$ in this case, which leads to $\psi_{r,q}\left(  \alpha^{k}%
\beta^{l}\right)  =0$. If $k+l=2N$ for some integer $N$,
(\ref{Divisibility condition psi red}) has solutions for $d$ precisely when
$\left\vert d\right\vert \leq R_{k,l}$ so that (\ref{Divisibility condition psi red}) reduces to $s+t=N+dw$ where $w=\frac
{1}{2}\left(  q-1\right)  $.

Consequently, we can simplify the right-hand side of (\ref{A form of psi red part 1}) to%
\begin{equation}
\psi_{r,q}\left(  \alpha^{k}\beta^{l}\right)  =\frac{\left(  -1\right)  ^{l}%
}{2^{2N}}\sum_{\left\vert d\right\vert \leq R_{k,l}}\sum_{s=0}^{k}\left(
-1\right)  ^{N+dw-s}\binom{k}{s}\binom{l}{N+dw-s}1_{q}.
\label{A form of psi red part 2}%
\end{equation}
For each $d$,%
\[
\sum_{s=0}^{k}\left(  -1\right)  ^{N+dw-s}\binom{k}{s}\binom{l}{N+dw-s}%
=\left[  \alpha^{N+dw}\right]  \left(  1+\alpha\right)  ^{k}\left(
1-\alpha\right)  ^{l}%
\]
so (\ref{A form of psi red part 2}) becomes
\begin{equation}
\psi_{r,q}\left(  \alpha^{k}\beta^{l}\right)  =\frac{\left(  -1\right)  ^{l}%
}{2^{2N}}\sum_{\left\vert d\right\vert \leq R_{k,l}}\left[  \alpha
^{N+dw}\right]  \left(  1+\alpha\right)  ^{k}\left(  1-\alpha\right)
^{l}1_{q}. \label{A form of psi red part 3}%
\end{equation}
Using Lemma \ref{Lemma palindromic property Krawtchouk},
\[
\left[  \alpha^{N+dw}\right]  \left(  1+\alpha\right)  ^{k}\left(
1-\alpha\right)  ^{l}=\left(  -1\right)  ^{k}\left[  \alpha^{N-dw}\right]
\left(  1+\alpha\right)  ^{k}\left(  1-\alpha\right)  ^{l}%
\]
for all $\left\vert d\right\vert \leq R_{k,l}$ so if $k$ is odd,
(\ref{A form of psi red part 3}) simplifies to $\psi_{r,q}\left(  \alpha
^{k}\beta^{l}\right)  =0.$

We may therefore write $k=2m$ and $l=2n$. In this case, $N=m+n$ and (\ref{A form of psi red part 3}) becomes%
\begin{align*}
\psi_{r,q}\left(  \alpha^{k}\beta^{l}\right)   &  =\frac{1}{2^{2m+2n}}%
\sum_{\left\vert d\right\vert \leq R_{k,l}}\left[  \alpha^{m+n+dw}\right]
\left(  1+\alpha\right)  ^{2m}\left(  1-\alpha\right)  ^{2n}1_{q}\\
&  =\left(  \sum_{\left\vert d\right\vert \leq R_{k,l}}\left[  \alpha
^{m+n+dw}\right]  \right)  \pi_{2m,2n}\operatorname{mod}p. \qedhere
\end{align*}
\end{proof}

Here is the special case that brings in the circular super Catalan numbers.

\begin{theorem}
\label{Thm formula red CIF SC}In the red geometry over $\mathbb{F}_{q}$ where
$q=p^{r}$ for some prime $p,$ suppose $k,l\in\mathbb{N}$ such that $0\leq
k+l<q-1$. Then
\[
\psi_{r,q}\left(  \alpha^{k}\beta^{l}\right)  =\left(  -1\right)  ^{n}%
\Omega(m,n)\operatorname{mod}p.
\]

\end{theorem}

\begin{proof}
From Theorem \ref{Theorem formula red CIF} we deduce that if $0\leq k+l<q-1$
then $R_{k,l}=0$ and so%
\begin{align*}
\psi_{r,q}\left(  \alpha^{k}\beta^{l}\right)   &  =\left[  \alpha
^{m+n}\right]  \pi_{2m,2n}\operatorname{mod}p=\left(  -1\right)  ^{n}\left[
\alpha^{m+n}\right]  \left(  -1\right)  ^{n}\pi_{2m,2n}\operatorname{mod}p\\
&  =\left(  -1\right)  ^{n}\Omega\left(  m,n\right)  \operatorname{mod}p
\end{align*}
where the last equation follows from the restatement of the result of
Georgiadis, Munemasa, and Tanaka \cite{Georgiadis} given as Theorem
\ref{Thm algebraic derivation Munemasa identity}.
\end{proof}

Now suppose $-1=\imath^{2}$ for some element $\imath\in\mathbb{F}_{q}$. From
Proposition \ref{Prop cyclic group structure of the unit circles}, the red and
blue unit circles are isomorphic groups, and so we may obtain an explicit
formula for $\psi_{b,q}$ which is a small tweak from that for $\psi_{r,q}$.

\begin{theorem}
\label{Thm formula blue CIF -1 square}In the blue geometry over $\mathbb{F}%
_{q}$ when $-1$ is a square, define $w\equiv\frac{1}{2}\left\vert
S_{b,q}\right\vert =\frac{1}{2}\left(  q-1\right)  $ and for natural numbers
$k$ and $l$, let $R_{k,l}\equiv\left\lfloor \frac{k+l}{\left\vert
S_{b,q}\right\vert }\right\rfloor .$ Then
\[
\psi_{b,q}\left(  \alpha^{k}\beta^{l}\right)  =\left\{
\begin{array}
[c]{cl}%
\left(
{\displaystyle\sum\limits_{\left\vert d\right\vert \leq R_{k,l}}}
\left[  \alpha^{m+n+dw}\right]  \right)  \left(  -1\right)  ^{n}\pi
_{k,l}\operatorname{mod}p & \text{if }k=2m\text{, }l=2n\\
0 & \text{otherwise.}%
\end{array}
\right.
\]

\end{theorem}

\begin{proof}
Using the bijection $\varrho$ in the proof of Proposition
\ref{Prop cyclic group structure of the unit circles}, we determine that%
\[
\psi_{b,q}\left(  \alpha^{k}\beta^{l}\right)  =-\sum_{\left[  x,y\right]  \in
S_{b,q}}x^{k}y^{l}=-\left(  -\imath\right)  ^{l}\sum_{\left[  x,y\imath
\right]  \in S_{r,q}}x^{k}\left(  y\imath\right)  ^{l}.
\]
By renaming the point $\left[  x,y\imath\right]  \in S_{r}\left(
\mathbb{F}_{q}\right)  $ to $\left[  x,y\right]  $,%
\begin{equation}
\psi_{b,q}\left(  \alpha^{k}\beta^{l}\right)  =-\left(  -\imath\right)
^{l}\sum_{\left[  x,y\imath\right]  \in S_{r,q}}x^{k}\left(  y\imath\right)
^{l}=-\left(  -\imath\right)  ^{l}\sum_{\left[  x,y\right]  \in S_{r,q}}%
x^{k}y^{l}=\left(  -\imath\right)  ^{l}\psi_{r,q}\left(  \alpha^{k}\beta
^{l}\right)  . \nonumber%
\end{equation}
The result then follows from the formula of $\psi_{r,q}\left(  \alpha^{k}%
\beta^{l}\right)  $ given in Theorem \ref{Theorem formula red CIF} above. If
either $k$ or $l$ is odd, then $\psi_{b,q}\left(  \alpha^{k}\beta^{l}\right)
=0$, and if $k=2m$ and $l=2n$, then%
\[
\psi_{b,q}\left(  \alpha^{k}\beta^{l}\right)  =\left(  -\imath\right)
^{2n}\psi_{r,q}\left(  \alpha^{2m}\beta^{2n}\right)  =\left(
{\displaystyle\sum\limits_{\left\vert d\right\vert \leq R_{k,l}}}
\left[  \alpha^{m+n+dw}\right]  \right)  \left(  -1\right)  ^{n}\pi
_{k,l}\operatorname{mod}p.\qedhere
\]

\end{proof}

\begin{corollary}
If $-1$ is a square in $\mathbb{F}_{q}$, then $\psi_{b,q}$ and $\psi_{r,q}$
are related in the following manner: for all $k,l\in\mathbb{N}$,%
\[
\psi_{b,q}\left(  \alpha^{k}\beta^{l}\right)  =\left\{
\begin{array}
[c]{cl}%
\psi_{r,q}\left(  \alpha^{k}\beta^{l}\right)  & \text{if }l\not \equiv
2\operatorname{mod}4\\
-\psi_{r,q}\left(  \alpha^{k}\beta^{l}\right)  & \text{if }l\equiv
2\operatorname{mod}4.
\end{array}
\right.
\]

\end{corollary}

\subsection{An explicit formula for $\psi_{b,q}$ when $\left(  \frac{-1}%
{q}\right)  =-1$ via complex summations}

We conclude our analysis by giving an explicit formula for $\psi_{b,q}$ when
$\left(  \frac{-1}{q}\right)  =-1$. We may think of $\mathbb{F}_{q}$ as being
embedded in the quadratic extension $\mathbb{C}_{b,q}$, and so $\psi
_{b,q}\left(  \alpha^{k}\beta^{l}\right)  $ might be identified as
\begin{equation}
\psi_{b,q}\left(  \alpha^{k}\beta^{l}\right)  1_{\mathbb{D}}=\sum_{\left[
x,y\right]  \in S_{b,q}}x^{k}y^{l}1_{\mathbb{D}}.
\label{Formula for psi_b,q as being embedded in D}%
\end{equation}

From identity (\ref{Conjugation generic complex number}), if $z=x1_{\mathbb{D}%
}+yi\in\mathbb{C}_{b,q}$ then
\[
x1_{\mathbb{D}}=\frac{1}{2}\left(  z+z^{\ast}\right)  \quad\text{and\quad
}y1_{\mathbb{D}}=-\frac{1}{2}\left(  z-z^{\ast}\right)  i.
\]
Now any $z\in S_{b,q}$ satisfies $z^{-1}=z^{\ast}$, so for any such $z$ we may
write%
\[
x1_{\mathbb{D}}=\frac{1}{2}\left(  z+z^{-1}\right)  \quad\text{and\quad
}y1_{\mathbb{D}}=-\frac{1}{2}\left(  z-z^{-1}\right)  i.
\]
This allows us to find a formula for $\psi_{b,q}\left(  \alpha^{k}\beta
^{l}\right)  1_{\mathbb{D}}$ in
(\ref{Formula for psi_b,q as being embedded in D}) as an algebraic function in
$z$ instead of a function in $x$ and $y$ given by%
\begin{align}
\psi_{b,q}\left(  \alpha^{k}\beta^{l}\right)  1_{\mathbb{D}}  &
=\sum_{\left[  x,y\right]  \in S_{b,q}}x^{k}y^{l}1_{\mathbb{D}}\nonumber\\
&  =\sum_{\left[  x,y\right]  \in S_{b,q}}\left(  x1_{\mathbb{D}}\right)
^{k}\left(  y1_{\mathbb{D}}\right)  ^{l}\nonumber\\
&  =\sum_{z\in S_{b,q}}\left(  \frac{1}{2}\left(  z+z^{-1}\right)  \right)
^{k}\left(  -\frac{1}{2}\left(  z-z^{-1}\right)  i\right)  ^{l}.
\label{Formula for psi_b,q in terms of powers of z}%
\end{align}
Therefore, a formula for $\psi_{b,q}\left(  \alpha^{k}\beta^{l}\right)  $ can
be found in reference to summations of integer powers of complex numbers over
$S_{b,q}$. Note that this is akin to an integration of the form $\int_{S^{1}%
}e^{in\theta}~d\theta$ from classical Fourier analysis.

We are now ready to present a formula for $\psi_{b,q}\left(  \alpha^{k}%
\beta^{l}\right)  $ which is parallel to the other Fourier summation formulas
that we have already obtained.

\begin{theorem}
\label{Thm formula blue CIF -1 not square} In the blue geometry over
$\mathbb{F}_{q}$ when $-1$ is not a square, define $w\equiv\frac{1}%
{2}\left\vert S_{b,q}\right\vert =\frac{1}{2}\left(  q+1\right)  $ and for
natural numbers $k$ and $l$, let $R_{k,l}\equiv\left\lfloor \frac{k+l}{\left\vert
S_{b,q}\right\vert }\right\rfloor $. We have that%
\[
\psi_{b,q}\left(  \alpha^{k}\beta^{l}\right)  =\left\{
\begin{array}
[c]{cl}%
\left(
{\displaystyle\sum\limits_{\left\vert d\right\vert \leq R_{k,l}}}
\left[  \alpha^{m+n+dw}\right]  \right)  \left(  -1\right)  ^{n}\pi
_{k,l}\operatorname{mod}p & \text{if }k=2m\text{, }l=2n\\
0 & \text{otherwise.}%
\end{array}
\right.
\]

\end{theorem}

\begin{proof}
By the embedding $\psi_{b,q}\left(  \alpha^{k}\beta^{l}\right)  \mapsto
\psi_{b,q}\left(  \alpha^{k}\beta^{l}\right)  1_{\mathbb{D}}$, from (\ref{Formula for psi_b,q in terms of powers of z}) we have that
\begin{align}
\psi_{b,q}\left(  \alpha^{k}\beta^{l}\right)  1_{\mathbb{D}}  &
=\frac{\left(  -1\right)  ^{l}}{2^{k+l}}i^{l}\sum_{z\in S_{b,q}}\left(
z+z^{-1}\right)  ^{k}\left(  z-z^{-1}\right)  ^{l}\nonumber\\
&  =\frac{\left(  -1\right)  ^{l}}{2^{k+l}}i^{l}\sum_{z\in S_{b,q}}\left(
\sum_{s=0}^{k}\dbinom{k}{s}z^{2s-k}\right)  \left(  \sum_{t=0}^{l}\left(
-1\right)  ^{l-t}\dbinom{l}{t}z^{2t-l}\right)\nonumber\\
&  =\frac{1}{2^{k+l}}i^{l}\sum_{s=0}^{k}\sum_{t=0}^{l}\left(  -1\right)
^{t}\dbinom{k}{s}\dbinom{l}{t}\sum_{z\in S_{b,q}}z^{2\left(  s+t\right)
-\left(  k+l\right)  }. \label{blue geometry eqn1b}%
\end{align}
In the third equality above, we interchanged the summation order.

By Proposition \ref{Prop cyclic group structure of the unit circles},
$S_{b,q}$ is a cyclic group of order $q+1$, and so by Lemma
\ref{Thm cyclic group geom series}, it follows that
\[
\sum_{z\in S_{b,q}}z^{2\left(  s+t\right)  -\left(  k+l\right)  }=\left\{
\begin{array}
[c]{cl}%
1_{\mathbb{D}} & \text{if }q+1\text{ divides }2\left(  s+t\right)  -\left(
k+l\right) \\
0 & \text{otherwise}%
\end{array}
\right.
\]
and so the terms in (\ref{blue geometry eqn1b}) will vanish except for those
$s$ and $t$ for which the equation
\begin{equation}
2\left(  s+t\right)  =k+l+d\left(  q+1\right)  \label{blue geometry eqn1c}%
\end{equation}
has solutions for $d$. Similar to the proof of Theorem
\ref{Theorem formula red CIF}, the case $k+l$ is odd gives $\psi_{b,q}\left(
\alpha^{k}\beta^{l}\right)  =0$ and if $k+l=2N$ for some integer $N$, then (\ref{blue geometry eqn1c}) has solutions for $d$ precisely when
$\left\vert d\right\vert \leq R_{k,l}$.

Each $d$ where $\left\vert d\right\vert \leq R_{k,l}$ then contributes to the
specific equation $s+t=N+dw$ where $w=\frac{1}{2}\left\vert S_{b,q}\right\vert
$. We can simplify (\ref{blue geometry eqn1b}) to%
\begin{equation}
\psi_{b,q}\left(  \alpha^{k}\beta^{l}\right)  1_{\mathbb{D}}=\frac{1}{2^{2N}%
}i^{l}\sum_{\left\vert d\right\vert \leq R_{k,l}}\sum_{s=0}^{k}\left(
-1\right)  ^{N+dw-s}\dbinom{k}{s}\dbinom{l}{N+dw-s}1_{\mathbb{D}}.
\label{blue geometry eqn3}%
\end{equation}
The same analysis pertains as in the red geometry: we have that%
\[
\sum_{s=0}^{k}\left(  -1\right)  ^{N+dw-s}\dbinom{k}{s}\dbinom{l}%
{N+dw-s}=\left[  \alpha^{N+dw}\right]  \left(  \left(  1+\alpha\right)
^{k}\left(  1-\alpha\right)  ^{l}\right)
\]
so if $k$ and $l$ are both odd, then the right-hand side of
(\ref{blue geometry eqn3}) vanishes and if $k=2m$ and $l=2n$, then we may
rewrite (\ref{blue geometry eqn3}) as
\begin{align*}
\psi_{b,q}\left(  \alpha^{k}\beta^{l}\right)  1_{\mathbb{D}}  &  =\frac
{1}{2^{2\left(  m+n\right)  }}i^{2n}\sum_{\left\vert d\right\vert \leq
R_{k,l}}\left(  \left[  \alpha^{m+n+dw}\right]  \left(  \left(  1+\alpha
\right)  ^{2m}\left(  1-\alpha\right)  ^{2n}\right)  \right)  1_{\mathbb{D}}\\
&  =\left(  \sum_{d\leq R_{k,l}}\left[  \alpha^{m+n+dw}\right]  \left(
-1\right)  ^{n}\left(  \left(  \frac{1+\alpha}{2}\right)  ^{2m}\left(
\frac{1-\alpha}{2}\right)  ^{2n}\right)  \right)  1_{\mathbb{D}}.
\end{align*}

It then follows that
\begin{align*}
\psi_{b,q}\left(  \alpha^{k}\beta^{l}\right)   &  =\sum_{d\leq R_{k,l}}\left[
\alpha^{m+n+dw}\right]  \left(  -1\right)  ^{n}\left(  \left(  \frac{1+\alpha
}{2}\right)  ^{2m}\left(  \frac{1-\alpha}{2}\right)  ^{2n}\right)  1_{q}\\
&  =\left(
{\displaystyle\sum\limits_{\left\vert d\right\vert \leq R_{k,l}}}
\left[  \alpha^{m+n+dw}\right]  \right)  \left(  -1\right)  ^{n}\pi
_{2m,2n}\operatorname{mod}p.\qedhere
\end{align*}
\end{proof}

\subsection{A combined formula for $\psi_{b,q}$}

Both formulas for $\psi_{b,q}$ are exactly the same except for the geometric
constant $w$ and consequently $R_{k,l}$. However, we can combine both formulas
into one to make it similar to the formula for $\psi_{g,q}$ and $\psi_{r,q}$
established in Theorem \ref{Thm CIF green group structure} and
\ref{Theorem formula red CIF}.

\begin{theorem}
[A unified formula for $\psi_{b,q}$]\label{Thm formula blue CIF combined}In
the blue geometry over a general finite field $\mathbb{F}_{q}$ of size
$q=p^{r}$ for some prime $p>2$ and positive integer $r$, define a geometric
constant
\[
w\equiv\frac{1}{2}\left\vert S_{b,q}\right\vert =\frac{1}{2}\left(  q-\left(
\frac{-1}{q}\right)  \right)  .
\]
For any natural numbers $k$ and $l$, let $R_{k,l}\equiv\left\lfloor \frac
{k+l}{\left\vert S_{b,q}\right\vert }\right\rfloor $. We have that%
\begin{equation}
\psi_{b,q}(\alpha^{k}\beta^{l})=\left\{
\begin{array}
[c]{cl}%
\left(
{\displaystyle\sum\limits_{\left\vert d\right\vert \leq R_{k,l}}}
\left[  \alpha^{m+n+dw}\right]  \right)  \left(  -1\right)  ^{n}\pi
_{k,l}\operatorname{mod}p & \text{if }k=2m\text{, }l=2n\\
0 & \text{otherwise.}%
\end{array}
\right.\nonumber
\end{equation}

\end{theorem}

Note that this formula is strikingly similar to the formula of $\psi_{r,q}$,
with the only difference being the presence of the factor $\left(  -1\right)
^{n}$. Moreover, we have the same analogue as in Theorem
\ref{Thm formula red CIF SC}.

\begin{theorem}
\label{Thm formula blue CIF SC}In the blue geometry over a field
$\mathbb{F}_{q}$, for any $k,l\in\mathbb{N}$ such that $0\leq k+l<q-1$, the
blue Fourier summation functional has the form
\[
\psi_{b,q}(\alpha^{k}\beta^{l})=\left\{
\begin{array}
[c]{cl}%
\Omega\left(  m,n\right)  \operatorname{mod}p & \text{if }k=2m\text{ and
}l=2n\\
0 & \text{otherwise.}%
\end{array}
\right.
\]

\end{theorem}

\begin{proof}
If either $k$ or $l$ is odd, then $\psi_{b,q}\left(  \alpha^{k}\beta
^{l}\right)  =0$ regardless of whether $-1$ is a square or not in
$\mathbb{F}_{q}$. Moreover, if $0\leq k+l<q-1$ then by Theorem
\ref{Thm formula blue CIF combined}, $R_{k,l}=\left\lfloor \frac
{k+l}{\left\vert S_{b,q}\right\vert }\right\rfloor =0$ and thus if $k=2m$ and
$l=2n$, by Theorem \ref{Thm algebraic derivation Munemasa identity} we consequently have
\[
\psi_{b,q}(\alpha^{k}\beta^{l})=\left[  \alpha^{m+n}\right]  \left(
-1\right)  ^{n}\pi_{2m,2n}=\Omega\left(  m,n\right)  \operatorname{mod}p. \qedhere
\]
\end{proof}

We conclude our analysis of the derivation of the formula for the Fourier
summation functional $\psi_{b,q}$ by highlighting the similarity of the
formula for $\psi_{b,q}\left(  \alpha^{k}\beta^{l}\right)  $ when $0\leq
k+l<q-1$ above with the characteristic $0$ integration formula by Baker and
Folland. We restate the result below for convenience, by renaming the
polynumber to $\alpha^{k}\beta^{l}$ to match with our notation:%
\[
\int_{S_2}\alpha^{k}\beta^{l}~d\mu=\left\{
\begin{array}
[c]{cl}%
\Omega\left(  m,n\right)  & \text{if }k=2m\text{ and }l=2n\\
0 & \text{otherwise.}%
\end{array}
\right.
\]

A pleasant feature of the analysis that we have presented here lies in its
flexibility to extend to general fields of characteristic zero. In a further
paper, we will show that the traditional unit circle integration theory over
the field of real numbers, that utilizes transcendental functions and infinite
processes and limits, can be replaced with a purely algebraic approach in
terms of a circular integral functional, and the pleasant connection to the
super Catalan numbers suggests that there is something combinatorial about
this theory of integration. With this algebraic approach, the connection
between algebraic combinatorics and harmonic analysis can be better
appreciated, and the characteristic zero and finite field situations can be
viewed more uniformly.

\section{Some examples of the computation of $\psi_{c,q}$}

Here are some examples of computing $\psi_{c,q}.$

\begin{example}
\label{Example q=17 revisited}Let $q=17$. From Example \ref{Example q = 17},
we know that $\psi_{r,17}\left(  \alpha^{6}\beta^{4}\right)  =3$. Now, using
the result in Theorem \ref{Theorem formula red CIF}, we confirm this:
\[
\psi_{r,17}\left(  \alpha^{6}\beta^{4}\right)  =\left(  -1\right)  ^{2}%
\Omega(3,2)\operatorname{mod}17=\frac{3}{256}\operatorname{mod}17=3.
\]
Now if $k=40$ and $l=24$, then $R_{k,l}=\left\lfloor \frac{k+l}{q-1}%
\right\rfloor =4$ and therefore%
\begin{align*}
\psi_{r,17}\left(  \alpha^{40}\beta^{24}\right)   &  =\sum_{\left\vert
d\right\vert \leq4}[\alpha^{32+8d}]\pi_{40,24}\operatorname{mod}17\\
&  =\left(  [\alpha^{0}]+[\alpha^{8}]+[\alpha^{16}]+[\alpha^{24}]+[\alpha
^{32}]+[\alpha^{40}]+[\alpha^{48}]+[\alpha^{56}]+[\alpha^{64}]\right)
\pi_{40,24}.
\end{align*}
Through inspection,%
\begingroup
\allowdisplaybreaks
\begin{align*}
\lbrack\alpha^{0}]\pi_{40,24}  &= [\alpha^{64}]\pi_{40,24}=\frac{1}{2^{64}} \\
\lbrack\alpha^{8}]\pi_{40,24}  &= [\alpha^{56}]\pi_{40,24}=\frac{90744}{2^{64}}\\
\lbrack\alpha^{16}]\pi_{40,24} &= [\alpha^{48}]\pi_{40,24}=-\frac{15\,156\,452}{2^{64}} \\
\lbrack\alpha^{24}]\pi_{40,24} &= [\alpha^{40}]\pi_{40,24}=-\frac{264\,053\,432}{2^{64}}\\
\lbrack\alpha^{32}]\pi_{40,24} &= \frac{1650\,887\,238}{2^{64}}%
\end{align*}
\endgroup
and thus
\begin{align*}
\psi_{r,17}\left(  \alpha^{40}\beta^{24}\right)   &  =\left(  2\left(
\frac{1}{2^{64}}+\frac{90744}{2^{64}}-\frac{15\,156\,452}{2^{64}}%
-\frac{264\,053\,432}{2^{64}}\right)  +\frac{1650\,887\,238}{2^{64}}\right)
\operatorname{mod}17\\
&  =\frac{33345}{2^{49}}\operatorname{mod}17=4.
\end{align*}

\end{example}

\begin{example}
\label{Example q=27 revisited}Let $q=27=3^{3}$. From Example
\ref{Example q = 27}, $\left(  \dfrac{-1}{3}\right)  =-1$ and $\psi
_{b,27}\left(  \alpha^{2}\right)  =2$ which agrees with Theorem
\ref{Thm formula blue CIF -1 not square}. Moreover, if $k=l=26$, then
$w=\frac{1}{2}\left(  q+1\right)  =14$ and $R_{k,l}=\left\lfloor \frac
{k+l}{q+1}\right\rfloor =1$. Therefore%
\begin{align*}
\psi_{b,27}\left(  \alpha^{26}\beta^{26}\right)   &  =\sum_{\left\vert
d\right\vert \leq1}[\alpha^{26+14d}]\left(  -1\right)  ^{13}\pi_{26,26}%
\operatorname{mod}3\\
&  =-\left(  [\alpha^{12}]+[\alpha^{26}]+[\alpha^{40}]\right)  \pi
_{26,26}\operatorname{mod}3\\
&  =-\left(  \frac{230\,230}{2^{52}}-\frac{10\,400\,600}{2^{52}}%
+\frac{230\,230}{2^{52}}\right)  \operatorname{mod}3=0.
\end{align*}

\end{example}

\section{Fourier Summation Program via Periodicity}

Since any polynumber $\pi\in\mathrm{Pol}_{2}\left(  \mathbb{F}_{q}\right)  $
admits a unique decomposition $\pi_{1}+\pi_{2}$ where $\pi_{1}\in
\mathcal{P}_{q}$, $\pi_{2}\in\ker\left(  \varepsilon\right)  $ by Proposition
\ref{Decomposition of Pol space} and $\psi_{c,q}\left(  \pi\right)
=\psi_{c,q}\left(  \pi_{1}\right)  $ from the proof of Theorem
\ref{Existence and Uniqueness of CIF on Pol space}, there ought to be a
simplification of the formulas for $\psi_{c,q}$ especially if the degree of
the input polynumber $\pi$ is large.

In this final section, we obtain a further simplification of the formulas for
$\psi_{r,q}$ and $\psi_{b,q}$ by reducing the degree of the input polynumber
$\pi$ so that it falls in the principal subspace $\mathcal{P}_{q}$.

\begin{theorem}
\label{Thm periodicity property}Let $\psi_{c,q}$ be the Fourier summation
functional in any of the three geometries. For all natural
numbers $k$, $l$, $r$, and $s$, we have that
\begin{equation}
\psi_{c,q}\left(  \alpha^{k+rq}\beta^{l+sq}\right)  =\psi_{c,q}\left(
\alpha^{k+r}\beta^{l+s}\right)  . \nonumber%
\end{equation}

\end{theorem}

\begin{proof}
The polynumber $\pi=\alpha^{k+rq}\beta^{l+sq}-\alpha^{k+r}\beta^{l+s}$
evaluates to the zero function since
\[
\varepsilon\left(  \pi\right)  \left(  x,y\right)  =x^{k+rq}y^{l+sq}%
-x^{k+r}y^{l+s}=x^{k}y^{l}\left(  \left(  x^{r}y^{s}\right)  ^{q}-x^{r}%
y^{s}\right)  =0.
\]
The conclusion follows by the Locality property.
\end{proof}

To illustrate the Periodicity property, we consider again the case $q=13$ in
the blue geometry and gather a table of values of $\psi_{b,13}\left(
\alpha^{2m}\beta^{2n}\right)  $ for various $m$ and $n$, since $\psi
_{b,13}\left(  \alpha^{k}\beta^{l}\right)  =0$ if either $k$ or $l$ is odd. By
Theorem \ref{Thm periodicity property}, all the information required to
determine $\psi_{b,13}\left(  \alpha^{2m}\beta^{2n}\right)  $ is contained in
the following finite list of values of $\psi_{b,13}\left(  \alpha^{2m}%
\beta^{2n}\right)  $ for $0\leq m,n\leq6$.
\[%
\begin{tabular}
[c]{c|ccccccc}%
$\psi_{b,13}\left(  \alpha^{2m}\beta^{2n}\right)  :m\backslash n$ & $0$ & $1$
& $2$ & $3$ & $4$ & $5$ & $6$\\\hline
$0$ & $1$ & $7$ & $2$ & $6$ & $2$ & $7$ & $3$\\
$1$ & $7$ & $5$ & $9$ & $4$ & $8$ & $4$ & $9$\\
$2$ & $2$ & $9$ & $5$ & $9$ & $4$ & $8$ & $4$\\
$3$ & $6$ & $4$ & $9$ & $5$ & $9$ & $4$ & $8$\\
$4$ & $2$ & $8$ & $4$ & $9$ & $5$ & $9$ & $4$\\
$5$ & $7$ & $4$ & $8$ & $4$ & $9$ & $5$ & $9$\\
$6$ & $3$ & $9$ & $4$ & $8$ & $4$ & $9$ & $5$%
\end{tabular}
\
\]

\begin{center}
\textbf{Table 4}: All the information to determine $\psi_{b,13}\left(
\alpha^{2m}\beta^{2n}\right)  $.
\end{center}

For example, $\psi_{b,13}\left(  \alpha^{100}\beta^{200}\right)  =\psi
_{b,13}\left(  \alpha^{16}\beta^{20}\right)  =\psi_{b,13}\left(  \alpha
^{4}\beta^{8}\right)  =4$, where we apply Theorem
\ref{Thm periodicity property} above twice: the first one with $\left(
k,l,r,s\right)  =\left(  9,5,7,15\right)  $ and the second one with $\left(
k,l,r,s\right)  =\left(  3,7,1,1\right)  $.

In general, in the red and blue geometries, it is sufficient to find all the
values of $\psi_{c,q}\left(  \alpha^{2m}\beta^{2n}\right)  $ where $0\leq
m,n\leq\frac{1}{2}\left(  q-1\right)  $ since $\psi_{c,q}\left(  \alpha
^{k}\beta^{l}\right)  =0$ if either $k$ or $l$ is odd by Theorem
\ref{Theorem formula red CIF} and \ref{Thm formula blue CIF combined}. We
first consider an explicit formula for $\psi_{r,q}$ and $\psi_{b,q}$ evaluated
at $\alpha^{2m}\beta^{2n}\in\mathcal{P}_{q}$ in the following proposition.

\begin{proposition}
\label{Thm periodicity red and blue -1 square}Suppose $\mathbb{F}_{q}$ is a
finite field with $q=p^{r}$ elements for some prime $p>2$ and positive integer
$r$.

\begin{enumerate}
\item In the red geometry, let $w=\frac{1}{2}\left(  q-1\right)  $. If $0\leq
m,n\leq w$, then%
\[
\psi_{r,q}\left(  \alpha^{2m}\beta^{2n}\right)  =\left(  \left(  -1\right)
^{n}\Omega\left(  m,n\right)  +2[\alpha^{m+n+w}]\pi_{2m,2n}+\frac{2}{4^{m+n}%
}\delta_{m+n,2w}\right)  \operatorname{mod}p.
\]

\item In the blue geometry where $\left(  \frac{-1}{q}\right)  =1$, let
$w=\frac{1}{2}\left(  q-1\right)  $. If $0\leq m,n\leq w$, then
\[
\psi_{b,q}\left(  \alpha^{2m}\beta^{2n}\right)  =\left(  \Omega\left(
m,n\right)  +2[\alpha^{m+n+w}]\left(  -1\right)  ^{n}\pi_{2m,2n}+\frac
{2}{4^{m+n}}\left(  -1\right)  ^{n}\delta_{m+n,2w}\right)  \operatorname{mod}%
p.
\]

\item In the blue geometry where $\left(  \frac{-1}{q}\right)  =-1$, let
$w=\frac{1}{2}\left(  q+1\right)  $. If $0\leq m,n\leq w-1$, then%
\[
\psi_{b,q}\left(  \alpha^{2m}\beta^{2n}\right)  =\left(  \Omega\left(
m,n\right)  +2[\alpha^{m+n+w}]\left(  -1\right)  ^{n}\pi_{2m,2n}\right)
\operatorname{mod}p.
\]

\end{enumerate}
\end{proposition}

\begin{proof}
In the red geometry, if $0\leq m+n<w$, it follows from Theorem
\ref{Theorem formula red CIF} that $R_{k,l}=\left\lfloor \frac{k+l}%
{q-1}\right\rfloor =0,$ so only the central term $\left(  -1\right)
^{n}\Omega\left(  m,n\right)  \operatorname{mod}p$ will appear in the formula.
If $w\leq m+n<2w$, then $R_{k,l}=1$ which gives us two additional terms
$\left[  \alpha^{m+n+w}\right]  \pi_{2m,2n}\operatorname{mod}p$ and $\left[
\alpha^{m+n-w}\right]  \pi_{2m,2n}\operatorname{mod}p$ in the formula of
$\psi_{r,q}$ and they are equal by Lemma
\ref{Lemma palindromic property Krawtchouk}. Finally, if $m+n=2w$, we have two
further additional terms which are equal, namely $\left[  \alpha
^{m+n+2w}\right]  \pi_{2m,2n}\operatorname{mod}p$ and $\left[  \alpha
^{m+n-2w}\right]  \pi_{2m,2n}\operatorname{mod}p$. Since
\[
\left[  \alpha^{m+n+2w}\right]  \pi_{2m,2n}\operatorname{mod}p=\left[
\alpha^{2m+2n}\right]  \left(  \frac{1+\alpha}{2}\right)  ^{2m}\left(
\frac{1-\alpha}{2}\right)  ^{2n}\operatorname{mod}p=\frac{1}{4^{m+n}}%
\delta_{m+n,2w}\operatorname{mod}p
\]
the conclusion follows.

The second identity follows by the result of Theorem
\ref{Thm formula blue CIF -1 square} that
\[
\psi_{b,q}\left(  \alpha^{2m}\beta^{2n}\right)  =\left(  -1\right)  ^{n}%
\psi_{r,q}\left(  \alpha^{2m}\beta^{2n}\right)
\]
whenever $-1$ is a square in $\mathbb{F}_{q}$. The third identity is done
similarly to the proof of the red geometry case.
\end{proof}

The following result is a general program for calculating $\psi_{r,q}\left(
\alpha^{2m}\beta^{2n}\right)  $ and $\psi_{b,q}\left(  \alpha^{2m}\beta
^{2n}\right)  $ that reduces the computational effort.

\begin{theorem}
[Fourier Summation Program]\label{Thm main simplification periodicity}Suppose
$\mathbb{F}_{q}$ is a finite field with $q=p^{r}$ elements for some prime
$p>2$ and positive integer $r$. Given any natural numbers $m$ and $n$, perform
the following procedure:

\begin{enumerate}
\item Set $m_{0}=2m$ and $n_{0}=2n$.

\item For each $k\in\mathbb{N}$, identify four integers $Q_{k}$, $R_{k}$,
$S_{k}$, and $T_{k}$ such that $Q_{k}$ and $S_{k}$ are the largest integers
satisfying%
\[
m_{k}=Q_{k}q+R_{k}\ \text{and\ }n_{k}=S_{k}q+T_{k}%
\]
for some natural numbers $R_{k}$ and $T_{k}$.

\item Set $m_{k+1}=Q_{k}+R_{k}$ and $n_{k+1}=S_{k}+T_{k}$. If $m_{k+1}\geq q$
or $n_{k+1}\geq q\,$, then repeat step 2. If not, then set $M=m_{k+1}$ and
$N=n_{k+1}$ and stop the procedure.
\end{enumerate}

We have that

\begin{enumerate}
\item The procedure will halt;

\item $M=2m^{\ast}$ and $N=2n^{\ast}$ for some $m^{\ast},n^{\ast}\in
\mathbb{N}$; and

\item In the red geometry,%
\[
\psi_{r,q}\left(  \alpha^{2m}\beta^{2n}\right)  =\left(  \left(  -1\right)
^{n^{\ast}}\Omega\left(  m^{\ast},n^{\ast}\right)  +2[\alpha^{m^{\ast}%
+n^{\ast}+\frac{q-1}{2}}]\pi_{2m^{\ast},2n^{\ast}}+\frac{2}{4^{m^{\ast
}+n^{\ast}}}\delta_{m^{\ast}+n^{\ast},q-1}\right)  \operatorname{mod}p.
\]
In the blue geometry, if $-1$ is a square in $\mathbb{F}_{q}$ then%
\[
\psi_{b,q}\left(  \alpha^{2m}\beta^{2n}\right)  =\left(  \Omega\left(
m^{\ast},n^{\ast}\right)  +2\left(  -1\right)  ^{n^{\ast}}[\alpha^{m^{\ast
}+n^{\ast}+\frac{q-1}{2}}]\pi_{2m^{\ast},2n^{\ast}}+\frac{2\left(  -1\right)
^{n^{\ast}}}{4^{m^{\ast}+n^{\ast}}}\delta_{m^{\ast}+n^{\ast},q-1}\right)
\operatorname{mod}p
\]
and if $-1$ is not a square in $\mathbb{F}_{q}$ then%
\[
\psi_{b,q}\left(  \alpha^{2m}\beta^{2n}\right)  =\left(  \Omega\left(
m^{\ast},n^{\ast}\right)  +2[\alpha^{m^{\ast}+n^{\ast}+\frac{q+1}{2}}]\left(
-1\right)  ^{n^{\ast}}\pi_{2m^{\ast},2n^{\ast}}\right)  \operatorname{mod}p.
\]

\end{enumerate}
\end{theorem}

\begin{proof}
If the procedure does not halt, then the sequence $\left(  m_{k}\right)  $ is
a strictly decreasing sequence of natural numbers since
\[
m_{k}=Q_{k}q+R_{k}>Q_{k}+R_{k}=m_{k+1}%
\]
by construction, which is impossible.

Now, since $m_{0}=Q_{0}q+R_{0}$ is even and $q$ is odd, it follows that
$Q_{0}$ and $R_{0}$ must both be odd or even, so $m_{1}=Q_{0}+R_{0}$ is even.
By a similar argument, $n_{1}=S_{0}+T_{0}$ is even. We can keep on going to
show that $m_{k}$ and $n_{k}$ are even. Since the procedure halts, $M$ and $N$
are both even.

Finally, by construction and Theorem \ref{Thm periodicity property},
\[
\psi_{c,q}\left(  \alpha^{m_{k}}\beta^{n_{k}}\right)  =\psi_{c,q}\left(
\alpha^{Q_{k}q+R_{k}}\beta^{S_{k}q+T_{k}}\right)  =\psi_{c,q}\left(
\alpha^{Q_{k}+R_{k}}\beta^{S_{k}+T_{k}}\right)  =\psi_{c,q}\left(
\alpha^{m_{k+1}}\beta^{n_{k+1}}\right)
\]
which shows that%
\[
\psi_{c,q}\left(  \alpha^{2m}\beta^{2n}\right)  =\psi_{c,q}\left(
\alpha^{m_{0}}\beta^{n_{0}}\right)  =\psi_{c,q}\left(  \alpha^{M}\beta
^{N}\right)  =\psi_{c,q}\left(  \alpha^{2m^{\ast}}\beta^{2n^{\ast}}\right)  .
\]
Since $0\leq m^{\ast},n^{\ast}\leq\frac{1}{2}\left(  q-1\right)  $, we can
apply Proposition \ref{Thm periodicity red and blue -1 square} to prove the
last assertion.
\end{proof}

We close this section with an example.

\begin{example}
Suppose $q=13$ and we wish to compute $\psi_{b,13}\left(  \alpha^{1000}%
\beta^{600}\right)  $. If we were to use the formula for $\psi_{b,q}$ in
Theorem \ref{Thm formula blue CIF combined}, we need to find particular
coefficients of the circular polynumber $\pi_{1000,600}$ of degree $1600$.
This is of course computationally expensive. However, with the aid of Theorem
\ref{Thm main simplification periodicity}, we can perform the Fourier
Summation Program applied to $m_{0}=1000$, $n_{0}=600$, and $q=13$ to get%
\[
\left(  m_{1},n_{1}\right)  =\left(  88,48\right)  \quad\left(  m_{2}%
,n_{2}\right)  =\left(  16,12\right)  \quad\left(  m_{3},n_{3}\right)
=\left(  4,12\right).
\]
Since now $m_{3},n_{3}<q$, the procedure halts. Therefore%
\begin{align*}
\psi_{b,13}\left(  \alpha^{1000}\beta^{600}\right)   &  =\psi_{b,13}\left(
\alpha^{4}\beta^{12}\right) \\
&  =\left(  \Omega\left(  2,6\right)  +2\left(  -1\right)  ^{6}\left[
\alpha^{14}\right]  \left(  \frac{1+\alpha}{2}\right)  ^{4}\left(
\frac{1-\alpha}{2}\right)  ^{12}\right)  \operatorname{mod}13\\
&  =9.
\end{align*}

\end{example}

\end{document}